\documentclass[12pt, a4, twoside]{amsart}

\usepackage[T1]{fontenc}
\usepackage{amsfonts}
\usepackage{fancyhdr}
\usepackage{amsmath,amssymb}
\usepackage{eufrak}
\usepackage[latin1]{inputenc}
\usepackage{latexsym}
\usepackage{amsbsy}
\usepackage{dsfont}
\usepackage[english,frenchb]{babel}
\usepackage[all,dvips]{xy}
\usepackage{graphicx}
\usepackage{amscd}
\usepackage{url}
\usepackage{vmargin}

\setmarginsrb{2.5cm}{2.5cm}{2.7cm}{2cm}{0cm}{1cm}{0cm}{1cm}

\newtheorem{theoreme}{Th\'{e}or\`{e}me}[section]

\newtheorem{definition}[theoreme]{D\'{e}finition}

\newtheorem{exemple}[theoreme]{Exemple}
\newtheorem{contre-exemple}[theoreme]{Contre-exemple}
\newtheorem{definitions}[theoreme]{D\'{e}finitions}

\newtheorem{lemme}[theoreme]{Lemme}

\newtheorem{proposition}[theoreme]{Proposition}

\newtheorem{remarque}[theoreme]{Remarque}

\numberwithin{equation}{section}
\newtheorem{theodinvol}[theoreme]{Th\'eor\`eme d'involutivit\'e générique}
\newtheorem{theodeCK}[theoreme]{Th\'eor\`eme de Cartan-K\"ahler}

\def\preuve{\smallskip\goodbreak{\it Preuve.~~--~\kern.3em}
     \ignorespaces}%
\def\qedbox{$\square$}%
\def\qed{\ifmmode\qedbox\else\unskip\ \hglue0mm\hfill
     \qedbox\smallskip\goodbreak\fi}%
\def\preuved#1{\smallskip\goodbreak{\it Preuve~{#1}.~~--~\kern.3em}
     \ignorespaces}%
\def\preuvedt#1{\smallskip\goodbreak{\it Preuve du th\'eor\`eme~\ref{#1}.~~--~\kern.3em}
     \ignorespaces}%
\def\preuvedts#1#2{\smallskip\goodbreak{\it Preuve du th\'eor\`eme~\ref{#1}~#2.~~--~\kern.3em}
     \ignorespaces}%
\def\idpreuve{\smallskip\goodbreak{\it Id\'ees de la preuve.~~--~\kern.3em}
     \ignorespaces}%
\def\idpreuvedt#1{\smallskip\goodbreak{\it Id\'ees de la preuve du th\'eor\`eme~\ref{#1}.~~--~\kern.3em}
     \ignorespaces}%
\newcommand{\CC}{\mathbb{C}}

\newcommand{\ZZ}{\mathbb{Z}}
\newcommand{\RR}{\mathbb{R}}

\newcommand{\D}{\mathcal{D}}

\newcommand{\E}{\mathcal{E}}
\newcommand{\J}{\mathcal{J}}
\newcommand{\M}{\mathcal{M}}
\newcommand{\I}{\mathcal{I}}

\newcommand{\K}{\mathcal{K}}

\newcommand{\G}{\mathcal{G}}
\newcommand{\F}{\mathcal{F}}

\renewcommand{\E}{\mathcal{E}}

\renewcommand{\O}{\mathcal{O}}

\newcommand{\dgl}{$\D$-groupo\"ide de Lie}
\newcommand{\dgls}{$\D$-groupo\"ides de Lie}

\begin{document}

\title[Feuilletages singuliers de codimension un]{Feuilletages singuliers de codimension un,\\  Groupo\"ide de Galois et int\'egrales premi\`eres}

\author{Guy Casale}

\address{Laboratoire \'Emile Picard, umr 5580 ufr mig, Universit\'{e} Paul Sabatier, 118 Route de Narbonne, 31062 Toulouse cedex 4, France}

\email{casale@picard.ups-tlse.fr}

%\subjclass{12H05, 32S65, 37F75, 53C10, 34Mxx}
%\keywords{feuilletage, groupo\"ide de Galois, structure g\'eom\'etrique transverse}

\maketitle

\tableofcontents

%ùùùùùùùùùùùùùùùùùùùùùùùùùùùùùùùùùùùùùùùùùùùùùùùùùùùùùùùùùùùùùùùùùùùùùùùùùùùùùùùùùùùùùùùùùùùùùùùùùùùùùùùùùùùùùùùùùùùùùùùùùùùùùùùùùùùùùùùùùùùùùùùùùùùùùùùùùùù
%ùùùùùùùùùùùùùùùùùùùùùùùùùùùùùùùùùùùùùùùùùùùùùùùùùùùùùùùùùùùùùùùùùùùùùùùùùùùùùùùùùùùùùùùùùùùùùùùùùùùùùùùùùùùùùùùùùùùùùùùùùùùùùùùùùùùùùùùùùùùùùùùùùùùùùùùùùùù
%ùùùùùùùùùùùùùùùùùùùùùùùùùùùùùùùùùùùùùùùùùùùùùùùùùùùùùùùùùùùùùùùùùùùùùùùùùùùùùùùùùùùùùùùùùùùùùùùùùùùùùùùùùùùùùùùùùùùùùùùùùùùùùùùùùùùùùùùùùùùùùùùùùùùùùùùùùùù
%ù
%ù      INTRODUCTION
%ù
%ùùùùùùùùùùùùùùùùùùùùùùùùùùùùùùùùùùùùùùùùùùùùùùùùùùùùùùùùùùùùùùùùùùùùùùùùùùùùùùùùùùùùùùùùùùùùùùùùùùùùùùùùùùùùùùùùùùùùùùùùùùùùùùùùùùùùùùùùùùùùùùùùùùùùùùùùùùù
%ùùùùùùùùùùùùùùùùùùùùùùùùùùùùùùùùùùùùùùùùùùùùùùùùùùùùùùùùùùùùùùùùùùùùùùùùùùùùùùùùùùùùùùùùùùùùùùùùùùùùùùùùùùùùùùùùùùùùùùùùùùùùùùùùùùùùùùùùùùùùùùùùùùùùùùùùùùù
%ùùùùùùùùùùùùùùùùùùùùùùùùùùùùùùùùùùùùùùùùùùùùùùùùùùùùùùùùùùùùùùùùùùùùùùùùùùùùùùùùùùùùùùùùùùùùùùùùùùùùùùùùùùùùùùùùùùùùùùùùùùùùùùùùùùùùùùùùùùùùùùùùùùùùùùùùùùù

\section*{Introduction}

En 2001, B. Malgrange a propos\'e dans \textit{Le groupo\"ide de Galois d'un feuilletage} \rm  (\cite{malgrangegroupdegalois}) une mani\`ere de g\'en\'eraliser le groupe de 
Galois diff\'erentiel, d\'efini pour une \'equation diff\'erentielle lin\'eaire, aux feuilletages singuliers. Deux changements importants par rapport aux 
th\'eories pr\'ec\'edentes apparaissent. Premi\`erement, on perd la structure de groupe alg\'ebrique remplac\'ee par celle de $\mathcal{D}$-groupo\"ide de 
Lie. Deuxi\`emement, alors que le groupe de Galois agissait sur les variables d\'ependantes (les inconnues des \'equations diff\'erentielles), le groupo\"ide
 de Galois agit sur l'ensemble des variables d\'ependantes et ind\'ependantes (l'espace portant le feuilletage donn\'e par les \'equations diff\'erentielles).\\ Le groupo\"ide de Galois est d\'efini par le syst\`eme maximal d'\'equations aux d\'eriv\'ees partielles qui v\'erifie les conditions suivantes :
\begin{itemize}
 \item Les flots des champs de vecteurs tangents au feuilletage sont des solutions de ce syst\`eme.
 \item Les inclusions de la d\'efinition \ref{defDGL} sont v\'erifi\'ees. Ces derni\`eres signifient que l'identit\'e est solution du syst\`eme, que la 
 compos\'ee de deux solutions est une solution et que l'inverse d'une solution est encore une solution.
\end{itemize}
D'apr\`es \cite{malgrangegroupdegalois}, un syst\`eme d'\'equations v\'erifiant ce dernier point est appel\'e \dgl.

Des d\'efinitions analogues bien qu'impr\'ecises ont \'et\'e esquiss\'ees par J. Drach \cite{Dr} et E. Vessiot \cite{Ve1}, \cite{Ve2}. Elles \'etaient 
bas\'ees sur la notion de syst\`emes automorphes d'int\'egrales premi\`eres du feuilletage. Les preuves de l'existence de tels syst\`emes semblent, 
malheureusement, incompl\`etes. En basant sa d\'efinition sur les propri\'et\'es dynamiques du feuilletage, sans aucune r\'ef\'erence aux int\'egrales 
premi\`eres, et en utilisant son th\'eor\`eme d'involutivit\'e g\'en\'erique \cite{M4}, B. Malgrange contourne les probl\`emes de d\'efinition et r\'esoud 
les probl\`emes d'existence.

Dans cet article, nous \'etudierons le groupo\"ide de Galois d'un germe de feuilletage de codimension un. Nous rappellerons dans une premi\`ere partie les
 d\'efinitions locales de $\mathcal{D}$-groupo\"ide de Lie \cite{malgrangegroupdegalois} et les r\'esultats relatifs aux $\mathcal{D}$-groupo\"ide de Lie au-dessus d'un disque 
 de $\mathbb{C}$ (\cite{C2}, \cite{malgrangenotes}). Ces derniers nous permettrons, dans une deuxi\`eme partie, de pr\'eciser la nature du syst\`eme d'\'equations
  aux d\'eriv\'ees partielles d\'efinissant le groupo\"ide de Galois du feuilletage et de caract\'eriser sa taille par un nombre : le rang transverse du 
  groupo\"ide de Galois du feuilletage. Ce nombre appartient \`a $\{0,1,2,3,\infty \}$.

Dans la troisi\`eme partie, nous montrerons les liens entre le groupo\"ide de Galois et l'existence de structures transverses. Nous discuterons suivant le 
rang transverse du groupo\"ide de Galois l'existence de structures m\'eromorphes transverses euclidienne (rang transverse \'egale \`a un), affine (rang 
transverse \'egale \`a deux) ou projective (rang transverse \'egale \`a trois). Nous \'enoncerons le r\'esultat en terme de suites de Godbillon-Vey 
m\'eromorphes pour le feuilletage, associ\'ees aux structures transverses. Le th\'eor\`eme \ref{suitegv} donne l'\'egalit\'e du rang transverse du 
groupo\"ide de Galois avec la longueur minimale des suites de Godbillon-Vey m\'eromorphes du feuilletage dans le cas des suites de longueur $1$, $2$ ou $3$.
 Ces r\'esultats ont \'et\'e annonc\'es dans \cite{malgrangepapierchinois} o\`u B. Malgrange prouve une des deux in\'egalit\'es d'une mani\`ere plus g\'eom\'etrique mais 
 essentiellement analogue \`a la notre. La preuve que nous donnons de l'autre in\'egalit\'e semble \^etre diff\'erente de celle de B. Malgrange.

Dans \cite{C2} nous caract\'erisons les germes de diff\'eomorphismes de $(\CC,0)$ solutions d'un $\mathcal{D}$-groupo\"ide de Lie au-dessus de 
$(\mathbb{C},0)$. Les germes de feuilletages de $(\mathbb{C}^2,0)$ \`a singularit\'es r\'eduites \'etant compl\`etement d\'ecrit par leurs holonomies, nous 
\'etudierons plus pariculi\`erement ces feuilletage dans la quatri\`eme partie. Nous expliquerons comment les r\'esultats de la partie pr\'ec\'edente 
compl\'et\'es par ceux de \cite{C2} redonne la caract\'erisation en terme d'invariants analytiques des germes de feuilletages de $(\mathbb{C}^2,0)$ \`a 
singularit\'es r\'eduites admettant une structure transverse m\'eromorphe affine ou projective. Nous retrouvons ainsi de mani\`ere ``galoisienne'' les 
r\'esultats de M. Berthier et F. Touzet \cite{BT} et ceux de F. Touzet \cite{T}.

Dans la cinqui\`eme partie, nous discuterons des diff\'erents types de transcendance d'int\'egrales premi\`eres (d\'efinition \ref{deftypedetransc}) du 
feuilletage caract\'eris\'es par le groupo\"ide de Galois (th\'eor\`eme \ref{typedetransc}). Suivant le rang transverse du groupo\"ide de Galois, le 
feuilletage admet une int\'egrale premi\`ere de type m\'eromorphe, Darboux, Liouville ou Riccati. Dans le cas d'un feuilletage dont le groupo\"ide de Galois
 est de rang transverse deux ou trois, le r\'esultat provient essentiellement d'un th\'eor\`eme de M. Singer \cite{S} et de sa version projective \cite{C} 
 via le th\'eor\`eme \ref{suitegv}. Pour les feuilletages de groupo\"ide de Galois de rang transverse un, la d\'emonstration consiste \`a prouver une version adapt\'ee
  du th\'eor\`eme de Singer. Le cas des feuilletages de groupo\"ide de Galois de rang transverse nul se traite d'un mani\`ere diff\'erente, nous construirons
   un quotient de l'espace des feuilles naturellement muni d'une structure de courbe analytique.

Dans la derni\`ere partie, nous d\'ecrirons quelques-unes des relations entre le groupo\"ide de Galois d'un feuilletage et la notion d'extension fortement 
normale de E. R. Kolchin (\cite{Ko}). Les int\'egrales premi\`eres construites dans la partie pr\'ec\'edente sont naturellement des \'el\'ements d'une 
extension fortement normale du corps des fonctions m\'eromorphes. Nous prouvons ensuite la r\'eciproque : si il existe une int\'egrale premi\`ere du 
feuilletage dans une extension fortement normale du corps des fonctions m\'eromorphes, le groupo\"ide de Galois est de rang transverse fini (th\'eor\`eme 
\ref{extfortnorm}).

\section{D\'{e}finitions et rappels}

Dans \cite{malgrangegroupdegalois} B. Malgrange d\'{e}finit la notion de $\mathcal{D}$-groupo\"{i}de de Lie et montre plusieurs propri\'{e}t\'{e}s de ces objets. Nous 
commen\c{c}ons par rappeler les d\'efinitions relatives aux \emph{D}-groupo\"ides de Lie au-dessus d'un polydisque $\Delta$ de $\mathbb{C}^n$.\\
L'espace des jets d'ordre $k$ d'applications inversibles de $\Delta$ dans $\Delta$ sera not\'e $J_k^*(\Delta)$. Le choix d'une coordonn\'ee $x$ sur $\Delta$
 permet de faire l'identification :
$$
J_k^*(\Delta) = \Delta\times\Delta\times GL_n(\mathbb{C}) \underset{2 \leq \vert\alpha\vert\leq k}{\times}\mathbb{C}^{n\vert\alpha\vert}
$$
avec les coordonn\'ees naturelles $(x_i,y_i,y_i^\alpha)$. Le multi-indice $\alpha $ appartient \`a $\mathbb{N}^n$ et on note $\vert\alpha\vert$ la somme de 
ses composantes. Nous noterons $\epsilon_j$ le multi-indice dont la seule coordonn\'ee non nulle est la $j$-i\`eme et est \'egale \`a $1$. On munit l'espace 
$\Delta \times \Delta$ du faisceau d'anneaux
$$
\mathcal{O}_{J_{k}^*(\Delta)}=\mathcal{O}_{\Delta\times\Delta}[y_i^{\alpha},\frac{1}{det(y^{\epsilon_j}_i)}]
$$
qui s'identifie \`a l'anneau des \'equations aux d\'eriv\'ees partielles, polynomiales en les d\'eriv\'ees, d'ordre inf\'erieur ou \'egal \`a $k$ ayant $n$
 variables ind\'ependantes et $n$ variables d\'ependantes. \\
\'Etant donn\'e un jet d'application $y(x)$ de $\Delta$ dans $\Delta$, nous regrouperons les d\'eriv\'ees $y_i^{\alpha}$ suivant leurs ordres
 $\vert\alpha\vert$. Nous noterons $y'$ la jacobienne de $y$ par rapport \`a $x$, $y''$ la hessienne (qui est \'el\'ement de 
 $S^2\mathbb{C}^n \otimes \mathbb{C}^n$) de $y$ par rapport \`a $x$, $y'''$ \'el\'ement de $ S^3\mathbb{C}^n \otimes \mathbb{C}^n $ la forme trilin\'eaire 
 des d\'eriv\'ees troisi\`emes \ldots \\
Les espaces $J_k^*(\Delta)$ sont de plus munis d'une structure de groupo\"ide par la donn\'ee
\begin{itemize}
\vspace{0.1cm}
\item de la projection source $s:J_k^*(\Delta)\rightarrow\Delta$ d\'efinie par $s(x,y,\ldots)=x$,
\vspace{0.1cm}
\item de la projection but $t:J_k^*(\Delta)\rightarrow\Delta$ d\'efinie par $t(x,y,\ldots)=y$,
\vspace{0.1cm}
\item d'une composition $c:J_k^*(\Delta)\times_{\Delta}J_k^*(\Delta)\rightarrow J_k^*(\Delta)$ d\'efinie sur les couples de jets $(h,g)$ tels que $t(h)=s(g)$
 par
$$
c((x,y,y',y'',\ldots),(y,z,z',z'',\ldots))=(x,z,z'y',z''(y',y') + z'y'', \ldots),
$$
\item d'une identit\'e, la sous-vari\'et\'e d\'efinie par les \'equations $x_i=y_i$ et $y^{\epsilon_j}_i=\delta^j_i$ pour $0\leq i,j \leq n$ et 
$y_i^{\alpha}=0$ pour $|\alpha| \geq 2$, donn\'ee par le plongement $e : \Delta \rightarrow J_k^*(\Delta)$ par $e(x)=(x,x,id,0,\ldots,0)$,
\vspace{0.1cm}
\item d'une inversion $i:J_k^*(\Delta)\rightarrow J_k^*(\Delta)$ qui \`a un jet $(x,y,y',y'',\ldots)$ fait correspondre le jet
$$
(y,x,(y')^{-1},-(y')^{-1}y''((y')^{-1},(y')^{-1}) ,\ldots),
$$
\end {itemize}
On a de plus $n$ d\'erivations $D_i:\mathcal{O}_{J_k^*(\Delta)}\rightarrow \mathcal{O}_{J_{k+1}^*(\Delta)}$ qui correspondent aux d\'erivations partielles
de fonctions compos\'ees. \'Etant donn\'ee une \'equation $E$ :
$$
D_i E=\frac{\partial E}{\partial x_i} + \sum_{\ell,\alpha} \frac{\partial E}{\partial y_\ell^\alpha}y_\ell^{ \alpha + \epsilon_i}.
$$
Toutes ces fl\`eches sont compatibles au projections naturelles $J_{k+1}^*(\Delta) \to J_{k}^*(\Delta)$ ce qui permet de les d\'efinir sur l'espace 
$\Delta \times \Delta$ muni de l'anneau $\O_{J^*(\Delta)}=\underset{\longrightarrow}{\text{lim}}\ \O_{J_{k}^*(\Delta)}$.
Les d\'efinitions suivantes sont issues de \cite{malgrangegroupdegalois}.
\begin{definition}
    \label{defnaive}
Un groupo\"ide d'ordre $k$ sur $\Delta$ est donn\'e par un id\'eal (= faisceau d'id\'eaux) coh\'erent $\mathcal{I}_k$ de $\mathcal{O}_{J_k^*(\Delta)}$ tel que :
  \begin{enumerate}
      \item $\mathcal{I}_k\subset Ker(e^*)$,
      \item $i^*\mathcal{I}_k \subset \mathcal{I}_k $,
      \item $c^*\mathcal{I}_k \subset \mathcal{I}_k\otimes_{\mathcal{O}_{\Delta}}1 + 1\otimes_{\mathcal{O}_{\Delta}}\mathcal{I}_k$ (la somme \'etant prise 
      comme somme d'id\'eaux).
  \end{enumerate}
\end{definition}
\noindent Cette d\'efinition est naturelle mais en pratique trop restrictive pour la d\'efinition de groupo\"ide de Galois telle qu'elle est donn\'ee dans 
la suite. Il faut alors utiliser la d\'efinition plus souple suivante.
\begin{definition}
    \label{defDGL}
Un $\D$-groupo\"ide de Lie sur $\Delta$ est donn\'e par un id\'eal r\'eduit $\mathcal{I}$ de $\mathcal{O}_{J^*(\Delta)}$ tel que
  \begin{itemize}
      \item tous les id\'eaux $\mathcal{I}_\ell=\mathcal{I}\cap \mathcal{O}_{J_\ell^*(\Delta)}$ sont coh\'erents,
      \item $\mathcal{I}$ soit stable par d\'erivation,
      \item pour tout ouvert relativement compact $U \subset \Delta$ il existe un entier $k$ et un ensemble analytique ferm\'e $Z$ dans $U$ tels que pour 
      tout $\ell\geq k$, $\widetilde{\mathcal{I}_\ell}=\mathcal{I}_\ell|_{U} $ v\'erifie \\
(i) les inclusions \emph{(1)} et \emph{(2)} de la d\'efinition \ref{defnaive},\\
(ii) l'inclusion \emph{(3)} de la d\'efinition \ref{defnaive} sur tout voisinage de $(x,y,z)\in (U-Z)\times (U-Z)\times (U-Z)$.
  \end{itemize}
\end{definition}

Dans cet article nous ne nous int\'eresserons qu'aux \dgls\ au-dessus d'un polydisque de taille arbitrairement petite. 
En nous pla\c{c}ant directement sur
 un polydisque plus petit que celui de d\'efinition, les points (\textit{i}) et (\textit{ii}) seront v\'erifi\'es sur tout le polydisque.

Une solution de $\mathcal{I}$ en $p \in \Delta$ est un morphisme $u$ de $\mathcal{O}_{J^*(\Delta)}/\mathcal{I}$ dans $\mathcal{O}_{\Delta,p}$ au-dessus de 
la restriction $\mathcal{O}_{\Delta} \to \mathcal{O}_{\Delta,p}$, tel que $u(D_iE)=\frac{\partial}{\partial x_i}u(E)$. En regardant $f=u(y)$, on obtient un 
germe $f:(\Delta,p) \to (\Delta,q)$ satisfaisant les \'equations diff\'erentielles de l'id\'eal $\mathcal{I}$. On d\'efinit de m\^eme les solutions formelles
 comme morphismes dans $\widehat{\mathcal{O}}_{\Delta,p}$. R\'eciproquement, un germe d'application inversible $f$, solution des \'equations diff\'erentielles
  engendrant $\mathcal{I}$ d\'efinit un morphisme $u$ par $u(y)=f$. Nous identifierons souvent, par abus de langage, un $\mathcal{D}$-groupo\"ide de Lie avec
   ses solutions formelles. En particulier nous dirons qu'un $\mathcal{D}$-groupo\"ide de Lie d'id\'eal $\I$ est inclus dans un second d'id\'eal $\J$ si 
   l'idéal $\I$ contient l'idéal $\J$ et qu'un \dgl\ contient un diff\'eomorphisme si ce dernier est solution des \'equations de $\I$.\\
Donnons quelques exemples de $\mathcal{D}$-groupo\"ides de Lie :

\begin{exemple}
\label{mero}
Le groupo\"ide d'invariance d'une fonction m\'eromorphe $\frac{P}{Q}$ est un
$\mathcal{D}$-groupo\"ide de Lie dont l'id\'eal est diff\'erentiablement
engendr\'e par
$$
Q(x)P(y)-Q(y)P(x).
$$
Les propri\'et\'es (1) et (2) d'un $\mathcal{D}$-groupo\"ide de Lie sont \'evidentes. La propri\'et\'e (3) est vraie en dehors du lieu d'ind\'etermination 
de $\frac{P}{Q}$. En effet, les \'egalit\'es
\begin{multline}
Q(z)[Q(x)P(y)-Q(y)P(x)] \nonumber \\
=Q(y)[Q(x)P(z)-Q(z)P(x)]-Q(x)[Q(z)P(y)-Q(y)P(z)] \nonumber
\end{multline}
et
\begin{multline}
P(z)[Q(x)P(y)-Q(y)P(x)]\nonumber \\
=P(y)[Q(x)P(z)-Q(z)P(x)]-P(x)[Q(z)P(y)-Q(y)P(z)] \nonumber
\end{multline}
donnent l'inclusion voulue tant que $P(z) \neq 0$ ou $Q(z) \neq 0$.

\end{exemple}

\begin{exemple}
Le groupo\"ide d'invariance d'un champ m\'eromorphe de
tenseurs, $T$, est un $\mathcal{D}$-groupo\"ide de Lie dont l'id\'eal est diff\'erentiablement engendr\'e par les composantes de $\Gamma^*T-T$. 
Les propri\'et\'es (1), (2) et (3) proviennent de l'\'egalit\'e
$$
\Gamma_2^*\Gamma_1^* T -T= \Gamma_2^*(\Gamma_1^* T - T) - (\Gamma_2^* T -T)
$$
L'ensemble $Z$ est alors inclus dans le lieu des p\^oles du champ $T$.
\end{exemple}

\begin{exemple}
Le groupo\"ide d'invariance d'un champ d'hyperplan donn\'e par une 1-forme
$\omega $ int\'egrable ($\omega \wedge d\omega = 0$) est un \dgl. Son id\'eal est
engendr\'e par les composantes de $\Gamma^* \omega \wedge \omega $ ou encore,
en coordonn\'ees dans lesquelles $\omega = \sum \omega_i dx_i$, par $\frac{(\Gamma^* \omega)_i}{(\Gamma^* \omega)_j}-\frac{\omega_i}{\omega_j}$. 
La troisi\`eme inclusion est v\'erifi\'ee en dehors du lieu d'annulation de $\omega$.

\end{exemple}

\noindent Les autres exemples que l'on pourrait donner sont des g\'en\'eralisations de ceux-ci en consid\'erant les groupo\"ides d'invariance 
(ou d'isom\'etries) de structures g\'eom\'etriques d'ordre sup\'erieur \`a un : voir \cite{Gr}, \cite{D}. \\

\'Etant donn\'e un syt\`eme d'\'equations aux d\'eriv\'ees partielles d'ordre $k$ : $\I_k$, l'id\'eal $pr_q\I_k$ des \'equations d'ordre $k+q$ obtenues par d\'erivations de $\I_k$
peut contenir des \'equations d'ordre $k$ n'appartenant pas à $\I_k$. Ceci nous interdit de consid\'erer les jets d'ordre $k$ solutions de $\I_k$
comme des jets de solutions formelles. Les syst\`emes diff\'erentiels ayant de bonnes propri\'et\'es d'int\'egrabilit\'e formelle 
($pr_q \I_k \cap \O_{J^*_{k+s}}=pr_{s}\I_k)$ sont les syst\`emes involutifs. Les th\'eor\`emes d'involutivité générique de Cartan-Kuranishi et 
de B. Malgrange nous assurent que n'importe quel syst\`eme différentiel est équivalent à un syst\`eme involutif en dehors d'une hypersurface de 
conditions initiales.

Soit $\I_k$ un système d'équations d'ordre $k$ tel que si $E \in \I_k \cap \O_{J_{k-1}(\Delta)}$ alors $D_iE \in \I_k$. On note $S_k$ la vari\'et\'e 
analytique d\'efinie par $\I_k$ d'anneau $\O_{S_k}= \O_{J_k}/\I_k$ et pour $E$ dans $\O_{J_k}$ on note $\delta E$ le symbole de $E$, c'est-\`a-dire sa 
diff\'erentielle modulo les $dx_i$ et les $dy_j^\alpha$ pour $|\alpha| \leq k-1$. Soient $E_1,\ldots E_p$ un syt\`eme d'\'equations qui engendre localement 
$\I_k$. Le premier prolongement de l'id\'eal est engendr\'e par les $E_\ell$ et les $D_i E_\ell$. Pour trouver un z\'ero de $pr_1 \I_k$ dans $J^*_{k+1}(\Delta)$ 
au-dessus d'un z\'ero de $\I_k$ dans $J^*_{k}(\Delta)$, il faut r\'esoudre un syst\`eme d'\'equations
$$
\sum_{j,|\alpha|=k}\frac{\partial E}{\partial y_j^{\alpha}} y_j^{\alpha+\epsilon_i}=\sum_{j,|\alpha|<k}\frac{\partial E}{\partial y_j^{\alpha}} y_j^{\alpha+\epsilon_i}.
$$
De m\^eme, pour que $pr_2 \I_k$ ait des z\'eros au-dessus de ceux de $pr_1 \I_k$, il faut r\'esoudre des \'equations de la forme 
$$
\sum_{j,|\alpha|=k}\frac{\partial E}{\partial y_j^{\alpha}} y_j^{\alpha+\epsilon_i+\epsilon_\ell}=*.
$$
La nature des prolongements successifs de $\I_k$ est donc control\'ee par les symboles.
 
On note $A[\xi]$ l'anneau des polyn\^omes en $\xi_1,\ldots \xi_n$ \`a coefficients dans un anneau $A$ et $A[\xi]_k$ l'espace des polyn\^omes homogènes de degré $k$.
\begin{definition}
Apr\`es la substitution de $\delta y_j^{\alpha}$ par $\xi^{\alpha}\delta y_j$, 
\begin{itemize} 
\item on appelle symbole d'ordre $k$ de $\I_k$ le $\O_{S_k}$-module $N_k$ engendr\'e dans\\ \hspace*{5mm}$\oplus_j \O_{S_k}[\xi]_k\delta y_j = \O_{S_k}[\xi]_k^m$ par les classes modulo $\I_k$ des $\delta f$ ;
\item on appelle symbole de $\I_k$ le $\O_{S_k}$-module gradué $N$ engendr\'e par \\ \hspace*{5mm} $N_k$ dans $\O_{S_k}[\xi]^m$ ;
\item on appelle module caract\'eristique de $\I_k$, le module gradué quotient \\ \hspace*{5mm} $M_k = \O_{S_k}[\xi]^m/ N$.
\end{itemize}  
\end{definition}

Dans la d\'efinition des systèmes différentiels ayant de bonnes propriétés de prolongement, les symboles interviendront par l'intermédiaire des \'evaluations ponctuelles du module caractéristique sur $S_k$. On demandera à ces derniers d'\^etre involutifs. 

\begin{definition}
Soit $M$ un $\CC[\xi]$-module gradu\'e. On dira que $M$ est $\ell$-involutif si il existe une base $(\eta_1,\ldots,\eta_n)$ de $\CC[\xi]_1$ v\'erifiant pour tout $q \geq \ell$ :\\
\indent Pour $i=1 \ldots n$, la multiplication par $\eta_i$
$$
M_q/(\eta_1,\ldots, \eta_{i-1})M_{q-1} \to M_{q+1}/(\eta_1,\ldots, \eta_{i-1})M_{q}
$$
\indent est injective et $M_q/(\eta_1,\ldots, \eta_{n})M_{q-1}=0$
\end{definition}

\begin{definition}
On dira qu'un système différentiel $\I_\ell$ est $\ell$-involutif si :
\begin{enumerate}
\item $S_\ell$ est lisse,
\item $M_\ell$ et $M_{\ell+1}$ sont localement libres,
\item en tout point $a$ de $S_\ell$, $M(a)$ est $\ell$-involutif,
\item $pr_1 S_{\ell} \to S_\ell$ est surjectif. 
\end{enumerate}
\end{definition}

Nous utiliserons les d\'efinitions précédentes uniquement \`a travers les trois th\'eor\`emes suivants. Pour leurs d\'emonstrations, nous renvoyons le lecteur à \cite{malgrangeprepubli}. 

\begin{theodeCK}
\label{theodeCK}
Si $\I_\ell$ est $\ell$-involutif alors $pr_1\I_\ell$ est $\ell+1$-involutif. Pour tout jet d'ordre $\ell$ solution de $\I_\ell$, il existe une solution convergente de $\I_\ell$ ayant ce jet d'ordre $\ell$.
\end{theodeCK}

\begin{theodinvol}[\cite{malgrangeprepubli}]
\label{theodinvol}
Soit $\I$ un id\'eal diff\'erentiel, r\'eduit de $\O_{J^*(\Delta)}$ tel que les id\'eaux $\I_k$ soient coh\'erents. Quitte \`a diminuer le polydisque $\Delta$, il existe un entier $\ell$
 et un sous-ensemble analytique ferm\'e de codimension un $Z_\ell \subset S_\ell$ v\'erifiant :
\begin{enumerate}
\item en dehors de $Z_\ell$, $\I$ est $\ell$-involutif,
\item pour tout entier $q$, en dehors de $Z_\ell$, $\I_{\ell+q}$ est le prolongement d'ordre $q$ de $\I_\ell$ et il n'y a aucune composante de $\I_{\ell+q}$ au-dessus de $Z_\ell$.   

\end{enumerate}
\end{theodinvol}

B. Malgrange montre parall\`element une version analytique du th\'eor\`eme de Ritt-Radenbush.
\begin{theoreme}
\label{malg-ritt-rad}
Soit $\I$ un id\'eal diff\'erentiel r\'eduit de $\O_{J^*(\Delta)}$. Quitte \`a r\'eduire $\Delta$, il existe un entier $\ell$ tel que l'id\'eal $\I$ soit l'id\'eal r\'eduit diff\'erentiablement engendr\'e par $\I_\ell$.
\end{theoreme}

Une premi\`ere cons\'equence de ces th\'eor\`emes est le th\'eor\`eme suivant.

\begin{theoreme}[\cite{malgrangegroupdegalois}]
\label{prolongementdesdgl}
Soit $\I_k$ un système différentiel d'ordre $k$. On note $\I$ l'id\'eal diff\'erentiel qu'il engendre et $\I^{\text{réd}}$ l'idéal réduit de $\I$. 
Supposons que $\I_k$ soit inclus dans $\ker e^*$, soit stable par $i^*$ et qu'il existe un ensemble analytique ferm\'e $Z$ de $\Delta$ tel
 que $\I_k$ d\'efinisse un groupo\"ide d'ordre $k$ en dehors de $Z$. Alors
\begin{enumerate}
\item $\I^{\text{réd}}$ est l'idéal d'un \dgl,
\item il existe un ensemble analytique ferm\'e $Z'$ de $\Delta$ tel que en dehors de $Z'$, $\I =\I^{\text{réd}}$.
\end{enumerate}
\end{theoreme}

\noindent Le th\'eorème de ``noeth\'erianité'' \ref{malg-ritt-rad} permet de montrer :

\begin{theoreme}[\cite{malgrangegroupdegalois}]
Soient $\I^{\alpha}$ des id\'eaux de \dgls. L'id\'eal réduit engendr\'e par la somme des $\I^{\alpha}$ est encore l'id\'eal d'un \dgl.
\end{theoreme}

Introduisons maintenant la notion de $\D$-alg\`ebre de Lie, \textit{i.e} la partie infinit\'esimal des \dgl. Commen\c{c}ons par
d\'efinir le crochet de Spencer sur les section de l'espace des jets d'ordre $k$ de champs de vecteurs $J_k(\Delta \to T \Delta)$
en suivant la construction diagonale \cite{K-S}.\\ 
Soit $R_k$ le fibr\'e sur $\Delta$ des jets d'ordre $k$ d'applications inversibles de $\Delta$ dans $(\CC^n,0)$. On a l'application suivante 
$$
\lambda :  R_k \times R_k \to J_k^*(\Delta) 
$$ 
definie par $ ( \varphi_2, \varphi_2) \mapsto \varphi_1 \circ \varphi_2^{-1} $. C'est le quotient de $R_k \times R_k$ sous l'action du groupe alg\'ebrique
des jets d'ordre $k$ de biholomorphisme de $(\CC^n,0)$ :
 $$ GL_n^{(k)}=J_k^*((\CC^n,0) \to (\CC^n,0))$$ par composition aux buts sur les deux facteurs.
Cette application induit une application du tangent vertical le long de la diagonale $ T (R_k \times R_k)/_{R_k} |_{diag}$ 
sur le tangent vertical le long de l'identit\'e $T (J_k^*(\Delta))/_{\Delta}|_{id}$ qui permet d'identifier les champs de 
vecteurs tangent \`a $J_k^*(\Delta)$ le long de l'identit\'e et verticaux aux champs de vecteurs tangents \`a 
$R_k \times R_k$ le long de la diagonale, verticaux et invariant sous l'action de $GL_n^{(k)}$.
On identifie ensuite le tangent vertical de $J_k^*(\Delta)$ le long de l'identit\'e à l'espace des jets d'ordre $k$ 
de champs de vecteurs $J_k(\Delta\rightarrow T\Delta)$ de la mani\`ere suivante. \`A un jet d'ordre $k$ de champ vertical en $a$, 
$\sum b_i^{\alpha}\frac{\partial}{\partial y_i^\alpha}$, on fait correspondre le jet $\sum b_i^\alpha (x-a)^\alpha \frac{\partial}{\partial x_i}$.
\begin{definition}
Gr\^ace aux construction ci-dessus, on peut identifier les sections de l'espace des jets 
$J_k(\Delta\rightarrow T\Delta)$ aux champs de vecteurs sur $R_k$ invariant sous l'action de $GL_n^{(k)}$.
 Ceci d\'efinit un crochet $[\ ,\ ]$ sur les sections de $J_k(\Delta\rightarrow T\Delta)$ en ramenant 
le crochet de Lie. C'est le crochet de Spencer. Ce crochet v\'erifie $[X^{(k)},Y^{(k)}]=[X,Y]^{(k)}$ pour les couples de champs de vecteurs sur $\Delta$ 
o\`u $(k)$ d\'esigne la section donn\'ee par le jet d'ordre $k$ d'un champ.
\end{definition}
\begin{definition}
\label{defalg}
La $\mathcal{D}$-alg\`ebre de Lie d'un $\mathcal{D}$-groupo\"ide de Lie d'id\'eal $\mathcal{I}$ est le lin\'earis\'e du $\mathcal{D}$-groupo\"ide le long de
 l'identit\'e. Elle est donn\'ee par l'id\'eal lin\'eaire et diff\'erentiel $\mathcal{L}(\mathcal{I})$ de $\mathcal{O}_{J^*(\Delta\rightarrow T\Delta)}$ 
 engendr\'e par les \'equations
$$
\mathcal{L}(E)=\sum_{i=1}^n \left( \frac{\partial E}{\partial y_i}(x,x,id,0,\ldots,0)a_i+\sum_{\vert \alpha \vert \leq k}\frac{\partial E}{\partial 
y_i^{\alpha}}(x,x,id,0,\ldots,0)a_i^{\alpha} \right)
$$
pour $E$ appartenant \`a $ \mathcal{I}_k$. Les $a_i$ sont les coordonn\'ees sur les fibres de $T\Delta$ induites par le choix de coordonn\'ees $x_i$ sur 
$\Delta$.
\end{definition}
Une solution de $\mathcal{L}(\mathcal{I})$ est un morphisme de $\mathcal{O}_{J^*(\Delta\rightarrow T\Delta)}/\mathcal{L}(\mathcal{I})$ dans $\mathcal{O}_{\Delta,p}^n$ 
commutant aux d\'erivations. Le choix de la coordonn\'ee $x$ sur $\Delta$ identifie $\mathcal{O}_{\Delta,p}^n$ aux germes de champs de vecteurs en $p$. 
Soient $(a_i(x))_{1 \leq i \leq n}$ les images des $(a_i)_{1\leq i\leq n}$ sous ce morphisme : le champ $\sum a_i(x)\frac{\partial}{\partial x_i}$ est 
appel\'e champ solution de $\mathcal{L}(\mathcal{I})$.
\begin{proposition}[\cite{malgrangegroupdegalois}]
Soit $\mathcal{I}$ l'id\'eal d'un $\mathcal{D}$-groupo\"ide de Lie. Pour tout entier $\ell$, le crochet de Spencer de deux sections solutions de 
$\mathcal{L}(\mathcal{I})_\ell = \mathcal{L}(\mathcal{I}) \cap \mathcal{O}_{J_\ell^*(\Delta\rightarrow T\Delta)} $ est une section solution.
\end{proposition}

\begin{definitions}
\text{ \\ }
\begin{enumerate}
   \item Une $\mathcal{D}$-alg\`ebre de Lie est donn\'ee par un id\'eal lin\'eaire diff\'erentiel $\mathcal{L}$ de $\mathcal{O}_{J^*(\Delta\rightarrow T\Delta)}$ 
   tel que pour tout $\ell$, les sections du fibr\'e vectoriel d\'efini par le lieu d'annulation de $\mathcal{L}_\ell$ dans 
   $J_{\ell}^*(\Delta \rightarrow T\Delta)$ soient stables par crochet de Spencer.
   \item Nous dirons qu'une $\mathcal{D}$-alg\`ebre de Lie est de rang $r$ lorsque le $\mathbb{C}$-espace vectoriel des champs solutions formelles en un point
    g\'en\'erique est de dimension $r$.
   \item Une $\mathcal{D}$-alg\`ebre de Lie sera dite int\'egrable si elle est la $\mathcal{D}$-alg\`ebre d'un $\mathcal{D}$-groupo\"ide de Lie.
   \item Un $\mathcal{D}$-groupo\"ide de Lie sera dit transitif lorsque les champs solutions de sa $\mathcal{D}$-alg\`ebre 
   de Lie \'evalu\'es en un point g\'en\'erique $x$ engendrent $T_x\Delta$ (ou, de mani\`ere \'equivalente, lorsque son id\'eal ne contient pas d'\'equation d'ordre z\'ero). 
\end{enumerate}
\end{definitions}

La classification et l'\'etude des $\mathcal{D}$-groupo\"ides de Lie au dessus d'un germe de disque de $\mathbb{C}$ est faite dans \cite{C2} \`a partir des 
notes de B. Malgrange \cite{malgrangenotes}. Nous en rappelons les r\'esultats principaux ci-dessous.

\begin{proposition}
\label{listealg}
Soit $\Delta$ un germe de disque de $\mathbb{C}$. Il y a exactement cinq types de $\mathcal{D}$-alg\`ebres de Lie sur $\Delta$, correspondant aux \'equations
 suivantes :
$$
\begin{array}{lcl}
\text{- rang 0: } a=0 & \text{not\'{e}e} & A_{0}\\
\text{- rang 1: } a'+\mu a=0  & \text{"} & A_{1}(\mu) \\
\text{- rang 2: } a''+\mu a'+\mu 'a=0  &\text{"}& A_{2}(\mu)  \\
\text{- rang 3: } a'''+\nu a'+\frac{\nu '}{2}a=0  & \text{"}& A_{3}(\nu) \\
\text{- rang $\infty$: \'equation nulle }&\text{"} & A_{\infty}\\
\end{array}
$$
les coefficents $\mu $ et $\nu $ \'{e}tant m\'{e}romorphes sur $\Delta $.
\end{proposition}

\begin{theoreme}
\label{integr1}
\text{    \\}
\begin{enumerate}
   \item L'\'equation $\nu(y){y'}^2+2\frac{y'''}{y'}-3\left(\frac{y''}{y'}\right)^2-\nu(x)=0$ d\'efinit l'unique $\mathcal{D}$-groupo\"{i}de de Lie ayant
    pour $\mathcal{D}$-alg\`ebre de Lie $ A_{3}(\nu )$ et sera not\'{e} $G_{3}(\nu)$. Dans ce cas, l'ensemble $Z$ (cf la d\'efinition \ref{defDGL}) est le 
    lieu des p\^oles de $\nu$.
   \item L'\'equation $\mu(y)y'+\frac{y''}{y'}-\mu(x)=0$ d\'efinit l'unique $\mathcal{D}$-groupo\"{i}de de Lie ayant pour $\mathcal{D}$-alg\`ebre de Lie 
   $ A_{2}(\mu )$ et sera not\'{e} $G_{2}(\mu )$. Dans ce cas, $Z$ est le lieu des p\^oles de $\mu$.
   \item La $\mathcal{D}$-alg\`ebre de Lie $A_{1}(\mu )$ n'est int\'{e}grable que lorsque $\mu $ a un p\^{o}le simple de r\'{e}sidu rationnel $\frac{p}{q}$.
    Les $\mathcal{D}$-groupo\"{i}des de Lie admettant $A_{1}(\mu )$ comme $\mathcal{D}$-alg\`ebre de Lie sont alors d\'{e}finis par les \'equations: 
   $$
   \gamma_k(y){y'}^{k}-\gamma_k(x)=0 \text{ avec } \gamma_k(x)=\exp
   (k \int \mu )(x)
   $$
   o\`u $k$ est un entier multiple de $q$. Ils seront not\'{e}s $G_{1}^{k}(\gamma_k )$. Dans ce cas, $Z$ est vide.
   \item Les $\mathcal{D}$-groupo\"ides de Lie ayant une $\mathcal{D}$-alg\`ebre de Lie de rang nul sont d\'efinis par une \'equation $h(x)-h(y)=0$ avec $h$
    holomorphe sur $\Delta$. Ici $Z$ est encore vide.
\end{enumerate}
\end{theoreme}

Nous rappelons maintenant les r\'esultats concernant les solutions $f:(\Delta,p)\to (\Delta,p)$ de ces $\mathcal{D}$-groupo\"ides de Lie.
\begin{proposition}
\label{flots}
Soit $a(x)\frac{d}{dx}$ un champ de vecteurs holomorphe sur un disque $\Delta$. 
Le germe de diff\'eomorphisme $f=exp(a(x)\frac{d}{dx})$ est solution des \dgls\
 suivants:
\begin{itemize}
\item $\mathcal{G}_1(-\frac{1}{a})$
\item $\mathcal{G}_2(\mu)$ avec $\mu$ solution m\'eromorphe de $a''+\mu a'+ \mu'a=0$
\item $\mathcal{G}_3(\nu)$ avec $\nu$ solution m\'eromorphe de $a'''+\nu a'+ \frac{\nu'}{2}a=0$
\end{itemize}
et seulement de ceux-ci.

\end{proposition}

\begin{exemple}
\label{exflot}
Les \dgls\ d'id\'eaux non nuls contenant le diff\'eomorphisme $\exp(x^2\frac{d}{dx})$ sont $G_1(-\frac{2}{x})$,
 $G_2(-\frac{c}{x^2}-\frac{2}{x})$ et
 $G_3(-\frac{c^2}{x^4})$ o\`u $c$ est une constante d'int\'egration.
\end{exemple}

%%%%%%%%%%%%%%%%%%%%%%%%%%%%%%%%%%%%%%%%%%%
Le cas des germes de diff\'eomorphismes formellement lin\'earisables est trait\'e par l'\'enonc\'e suivant.
\begin{theoreme}[\cite{malgrangenotes}]
\label{pasdepetitdiv}
Un diff\'eomorphisme $f$ formellement lin\'earisable est solution d'un $\mathcal{D}$-grou\-po\"ide de rang fini si et seulement si $f$ est analytiquement 
lin\'{e}arisable.
\end{theoreme}
\noindent Dans ce cas $f$ est solution d'un groupo\"{i}de de rang 1 : $G_1(\frac{h'}{h})$ o\`u $h$ est la lin\'earisante analytique de $f$.

%%%%%%%%%%%%%%%%%%%%%%%%%%%%%%%%%%%%%%%%%%%
Le cas d'un diff\'eomorphisme $f$ parabolique, c'est-\`a-dire vérifiant $f'(0)=1$, n\'ecessite la mise en place des notations suivantes (\cite{MR2}).
Nous noterons 
$$a_{k,\lambda}(x)=2i\pi\frac{x^{k+1}}{1+\lambda x^k} \text{ et } X_{k,\lambda}=a_{k,\lambda}\frac{d}{dx}.$$
 Les germes de diff\'eo\-morphismes 
$g_{k,\lambda}=\exp(X_{k,\lambda})$ sont les formes normales formelles des germes de diff\'eo\-morphismes paraboliques. Ils admettent pour int\'egrales 
premi\`eres $H_{k,\lambda}(x)=x^{-\lambda}e^{1/kx^k}.$ Pour tout $f(z)=z+cz^{k+1}+\ldots $, il existe un diff\'eo\-morphisme formel tangent \`a l'identit\'e
 $\widehat{h}$ qui conjugue $f$ \`a $g_{k,\lambda}$. Un th\'eor\`eme d'\'Ecalle, Voronin, Martinet-Ramis (\cite{E},\cite{malgrangebourbaki}) dit que 
 $\widehat{h}$ est k-sommable de somme $(U_i,h_i)$ sur $2k$ secteurs $U_i$. On note $inv(f)=(h_i(U_i)\cap h_{i+1}(U_{i+1}),h_{i+1}\circ h^{-1}_i)$ les 
 invariants analytiques d'\'Ecalle et Voronin d\'efinis \`a automorphisme pr\`es sur le disque du mod\`ele formel et $s(f)=(U_i \cap U_{i+1},h^{-1}_{i+1}
 \circ h_i)$ leurs analogues sur $\Delta$.

\begin{proposition}
\label{5}
\text{ \\}
\begin{enumerate}
\item Si $f$ est solution d'un $\mathcal{D}$-groupo\"ide de Lie, $s(f)$ est solution du m\^{e}me $\mathcal{D}$-groupo\"ide de Lie.
\item Soit $f$ de forme normale $g_{k,\lambda }$ et d'invariant analytique $inv(f)$. Le diff\'eo\-morphisme $f$ est solution d'un $\mathcal{D}$-groupo\"ide 
de rang $r$ si et seulement si il existe un $\mathcal{D}$-groupo\"ide de rang $r$ admettant $g_{k,\lambda}$ et $inv(f)$ comme solutions.
\end{enumerate}
\end{proposition}

Le th\'eor\`eme suivant donne la liste des germes de diff\'eomorphismes paraboliques solutions d'un \dgl\ de rang fini. Pour cela nous utiliserons les 
invariants ``g\'eo\-m\'etriques'' $(\varphi_{i,i+1})$ de Martinet-Ramis d\'efinis par $\varphi_{i,i+1} \circ (H_{k,\lambda} \circ h_{i+1})=H_{k,\lambda} 
\circ h_i$. Ils sont alternativement d\'efinis au voisinage de 0 et de l'infini sur les sph\`eres du chapelet d\'ecrivant l'espace des orbites de $f$ 
(\cite{MR2}).

\begin{theoreme}
\label{theoimp}
\text{ \\}
\begin{enumerate}
\item Un diff\'{e}omorphisme tangent \`{a} l'identit\'{e} est solution d'un $\mathcal{D}$-groupo\"{i}de de Lie de rang 3 si et seulement si il existe un entier 
positif $p$ tel que ses invariants g\'eom\'etriques soient de la forme $\frac{\tau}{\sqrt[p]{1+a_{i}\tau^{p}}}$ en $0$ et $\tau\sqrt[p]{1+\frac{b_{i}}{\tau^{p}}}$
 en $\infty $. Suivant la terminologie de J. \'Ecalle, de tels diff\'{e}omorphismes seront appel\'{e}s binaires.
\item Un diff\'{e}omorphisme tangent \`{a} l'identit\'{e} est solution d'un $\mathcal{D}$-groupo\"{i}de de Lie de rang 2 si et seulement si il existe un entier 
$p$ tel que ses invariants g\'eom\'etriques soient tous de la forme $\frac{\tau}{\sqrt[p]{1+a_{i}\tau^{p}}}$ en $0$ et $\tau$ en $\infty$ ou bien tous de la
 forme $\tau$ en $0$ et $\tau\sqrt[p]{1+\frac{b_{i}}{\tau^{p}}}$ en $\infty $. De tels diff\'{e}omorphismes seront appel\'{e}s unitaires.
\end{enumerate}
\end{theoreme}
Le cas des diff\'eomorphismes r\'esonnants, $f(x)=e^{2i\pi\frac{p}{q}}x + \ldots$, se d\'eduit de ce th\'eor\`eme en consid\'erant $f^{\circ q}$.

%%%%%%%%%%%%%%%%%%%%%%%%%%%%%%%%%%%%%%%%%%%%%%%%%%%%%%%

\section{Le groupo\"ide de Galois d'un feuilletage}

Soit $\mathcal{F}$ un germe de feuilletage holomorphe singulier de codimension $q$ d\'efini par des formes $\omega_1,\ldots \omega_q$ ind\'ependantes sur le
 corps des fonctions m\'eromorphes v\'erifiant les conditions d'int\'egrabilit\'e de Frobenius $d\omega_i \wedge \omega_1\wedge \ldots \wedge \omega_q = 0$.

Le groupo\"ide d'holonomie d'un feuilletage est un objet transverse d\'efini en dehors des singularit\'es du feuilletage et il ne refl\`ete pas compl\`etement
 la complexit\'e des singularit\'es. Le groupo\"ide de Galois d'un feuilletage est la ``cl\^oture de Zariski'' du groupo\"ide d'holonomie au sens de la 
 d\'efinition suivante :
\begin{definition}[\cite{malgrangegroupdegalois}]
Soit $\mathcal{F}$ un germe de feuilletage holomorphe singulier de $(\mathbb{C}^{n},0)$. Son groupo\"{i}de de Galois, $\mathcal{G}al(\mathcal{F})$, est le 
plus petit $\mathcal{D}$-groupo\"{i}de de Lie tel que les champs tangents \`{a} $\mathcal{F}$ soient solutions de sa $\mathcal{D}$-alg\`ebre de Lie.
\end{definition}
 L'existence de cet objet est une cons\'equence direct du th\'eor\`eme \ref{malg-ritt-rad}. Quelques unes de ses propri\'{e}t\'{e}s ont \'{e}t\'{e} 
 \'{e}tablies dans \cite{malgrangegroupdegalois}.
 
\begin{remarque}
Les \'equations $\omega_i(X)=0$ d\'efinissant les champs tangents au feuilletage sont lin\'eaires d'ordre z\'ero. L'hypoth\`ese
d'int\'egrabilit\'e de Frobenius se traduit par le fait qu'elles d\'efinissent un faisceau en sous-$\mathbb{C}$-alg\`ebre de Lie de $T\Delta$. L'id\'eal engendr\'e par
 ces \'equations et toutes leurs d\'eriv\'ees d\'efinit une $\mathcal{D}$-alg\`ebre de Lie, \emph{\cite{malgrangegroupdegalois}}. Dans la plupart des cas cette
  $\mathcal{D}$-alg\`ebre de Lie n'est pas int\'egrable.
\end{remarque}

Le $\mathcal{D}$-groupo\"ide de Lie $Aut(\mathcal{F}_{\omega})$ des automorphismes du feuilletage a pour solutions les germes de diff\'eomorphismes 
$\Gamma$ v\'erifiant $\Gamma^{\ast}\omega_i \wedge \omega_1 \wedge \ldots \wedge \omega_q =0$ . Il contient le groupo\"ide de Galois du feuilletage.

\begin{definition}
Un $\mathcal{D}$-groupo\"ide de Lie admissible pour $\mathcal{F}$ est un $\mathcal{D}$-groupo\"ide de Lie contenu dans $Aut(\mathcal{F}_{\omega})$ dont 
la $\mathcal{D}$-alg\`ebre de Lie contient la $\mathcal{D}$-alg\`ebre de Lie d\'efinie par le feuilletage.
\end{definition}

Un tel \dgl contient donc le groupo\"ide de Galois du feuilletage. 
Nous allons \'{e}tudier l'expression locale des \'equations d'un \dgl\ admissible sur un ouvert de redressement $U$. Nous noterons $t$ les coordonn\'ees 
transverses et $z$ les coordonn\'ees tangentes. Les $\mathcal{D}$-groupo\"ides de Lie au-dessus de $U$ sont donn\'es par des \'equations sur les 
diff\'erents espaces de jets
$$
J^*_k(U) = U \times U \times GL_n(\mathbb{C}) \underset{2 \leq \vert\alpha\vert\leq k}{\times}\mathbb{C}^{n\vert\alpha\vert}.
$$
Nous noterons $(t,z)$ les coordonn\'ees sur le premier $U$, $(T,Z)$ les m\^emes coordonn\'ees sur le second $U$ et $\frac{\partial^{|\alpha| + |\beta|}}
{\partial t^{\alpha} \partial z^{\beta}}T$, $\frac{\partial^{|\alpha| + |\beta|}}{\partial t^{\alpha} \partial z^{\beta}}Z$
 les coordonn\'ees naturellement induites  
sur l'espace des jets.

\begin{lemme}
\label{local}
\text{ \\ }
\begin{enumerate}
   \item Soient $\mathcal{F}$ un feuilletage sur $\Delta$ et $(U,(t,z))$ un ouvert muni de coordonn\'ees redressantes où $t$ désigne les coordonn\'ees 
   transverses et $z$ les coordonn\'ees tangentes. L'id\'eal d'un $\mathcal{D}$-groupo\"ide de Lie admissible pour $\mathcal{F}$ est engendr\'e sur $U$ 
   par des \'equations $$\frac{\partial T_j}{\partial z_i}=0 \text{ et } E_i(t,T,\ldots \frac{\partial^{|\alpha|} T}{\partial t^{\alpha} })$$ o\`u les $ E_i$ 
   sont les \'equations d'un $\mathcal{D}$-groupo\"ide de Lie au-dessus du polydisque transverse $t(U)$.
   \item Le rang du $\mathcal{D}$-groupo\"ide transverse ainsi obtenu est ind\'ependant de la carte choisie.
\end{enumerate}
\end{lemme}

\begin{definition}
Ce rang est appel\'e rang transverse du $\mathcal{D}$-groupo\"ide admissible.
\end{definition}
\begin{preuve}
Soit $\mathcal{I}$ l'id\'eal d'un $\mathcal{D}$-groupo\"ide de Lie admissible. Il s'obtient en compl\'etant les \'equations de $Aut(\mathcal{F}_{\omega})$ :
 $\frac{\partial T_j}{\partial z_i}=0$ par des \'equations suppl\'ementaires $F_i(t,z,T,Z,\ldots)$. Toute solution $\Gamma$ de $\mathcal{I}$ se factorise
  sous la forme
$\Gamma^{trans}\circ \Gamma^{tang}$ avec
$$
\Gamma^{trans}(t,z) = (T(t),z) \text{ et }
 \Gamma^{tang}(t,z) = (t,Z(t,z)).
$$
Comme $\mathcal{I}$ d\'efinit un $\mathcal{D}$-groupo\"ide admissible, toutes les transformations de la forme $\Gamma^{tang}$ sont solutions de $\mathcal{I}$.
 Pour toute solution $\Gamma$ de $\mathcal{I}$, $\Gamma^{trans}$ est aussi solution de $\mathcal{I}$ et v\'erifie de plus les \'equations de l'id\'eal 
 diff\'erentiel engendr\'e par $(Z-z)$. L'id\'eal diff\'erentiel $\mathcal{I}+(Z-z)$ d\'ecrit un \dgl\ et est engendr\'e par des \'equations de la forme
$$
\frac{\partial T_j}{\partial z_i},\ (Z_j-z_j),\ E_k(z,t,T,\ldots,\frac{ \partial^{|\alpha|} T }{ \partial t^{|\alpha|}}\ldots ).
$$
Les translations $\tau(t,z)=(t,z+a)$ sont solutions de $\mathcal{I}$. La conjugaison par $\tau$ laisse $\mathcal{I}$ invariant ainsi que $(Z-z)$. 
L'id\'eal engendr\'e par les \'equations ci-dessus est donc \'egal à l'id\'eal engendr\'e diff\'erentiablement par
$$
\frac{\partial T_j}{\partial z_i},\ (Z_j-z_j),\ E_k(z_0,t,T,\ldots,\frac{ \partial^{|\alpha|} T }{ \partial t^{|\alpha|}}\ldots ).
$$
 Ceci permet de choisir des g\'en\'erateurs $E_i$ ind\'ependants de $z$. Consid\'erons l'id\'eal $\mathcal{J} = \left( \frac{\partial T_j}{\partial z_i},
 E_k(t,T,\ldots,\frac{\partial ^{|\alpha|} T }{\partial t^{|\alpha|}})\ldots \right)$. Comme $\mathcal{J}+(Z-z)=\mathcal{I}+(Z-z)$ est l'id\'eal d'un \dgl,
  l'id\'eal engendr\'e au-dessus du disque tranverse par les \'equations $E_i(t,T,\ldots,\frac{d^k T }{dt^k}\ldots)$ d\'efinit un \dgl. Pour un entier 
  $\ell$ assez grand, un jet d'ordre $\ell$ de transformation $\Gamma=\Gamma^{trans}\circ \Gamma^{tang}$ est solution de $\mathcal{J}_\ell$ si et seulement 
  si $\Gamma^{trans}$ est solution de $\mathcal{I}_\ell$. La partie tangente $\Gamma^{tang}$ \'etant toujours solution de $\mathcal{I}_\ell$, 
  $\mathcal{I}_\ell$ et $\mathcal{J}_\ell$ ont m\^emes solutions parmi les jets d'ordre $\ell$. Comme ce sont des id\'eaux r\'eduits, ils sont \'egaux. 
  Les th\'eor\`emes \ref{theodinvol} et \ref{prolongementdesdgl} permettent de conclure que $\mathcal{I}=\mathcal{J}$.

\end{preuve}

Ce lemme ram\`ene l'\'etude locale en un point r\'egulier du feuilletage de son groupo\"ide de Galois à la compréhension des \dgls\ d\'efinis au-dessus du 
polydisque transverse et de leurs prolongements analytiques sur le polydisque de d\'efinition du feuilletage. En particulier pour les feuilletages de codimension un, 
le groupo\"ide de Galois sera d\'ecrit par des \dgls\ au-dessus d'un disque.

%ùùùùùùùùùùùùùùùùùùùùùùùùùùùùùùùùùùùùùùùùùùùùùùùùùùùùùùùùùùùùùùùùùùùùùùùùùùùùùùùùùùùùùùùùùùùùùùùùùùùùùùùùùùùùùùùùùùùùùùùùùùùùùùùùùùùùùùùùùùùùùùùùùùùùùùùùùùùùùùùùùùùùùùùùùùùùùùùùùùùùùùùùùùùùùùùùùùùùùùùùùùùùùùùùùùùùùùùùùùùùùùùùùùùùùùùùùùùùùùùùùùùùùùùùùùùùùùùùùùùùùùùùùùùùùùùùùùùùùùùùùùùùùùùùùùùùùùùùùùùùùùùùùùùùùùùùùùùùùùùùùùùùùùùùùùùùùùùùùùùùùùùùùùùùùùùùùùùùùùùùùùùùùùùùùùùùùùùùùùùùùùùùùùùùùùùùùùùùù
%ù          CHAPITRE 3
%ù
%ùùùùùùùùùùùùùùùùùùùùùùùùùùùùùùùùùùùùùùùùùùùùùùùùùùùùùùùùùùùùùùùùùùùùùùùùùùùùùùùùùùùùùùùùùùùùùùùùùùùùùùùùùùùùùùùùùùùùùùùùùùùùùùùùùùùùùùùùùùùùùùùùùùùùùùùùùùùùùùùùùùùùùùùùùùùùùùùùùùùùùùùùùùùùùùùùùùùùùùùùùùùùùùùùùùùùùùùùùùùùùùùùùùùùùùùùùùùùùùùùùùùùùùùùùùùùùùùùùùùùùùùùùùùùùùùùùùùùùùùùùùùùùùùùùùùùùùùùùù

\section{Groupo\"ides de Galois et Suites de Godbillon-Vey}

Dans ce paragraphe, $\mathcal{F}_\omega$ d\'esigne le feuilletage holomorphe
singulier donn\'e par une 1-forme $\omega$ int\'egrable ($\omega \wedge d\omega = 0$) sur un polydisque $\Delta$ dans $\mathbb{C}^n$. On pourra supposer
 que le lieu singulier de $\omega$ est de codimension deux.
 
\begin{definition}[\cite{GV}]
Une suite de Godbillon-Vey pour $\omega$ est une suite de 1-formes m\'{e}romorphes $\omega_{1},\omega_{2},\ldots,\omega_{n},\ldots$ telles que:

$$
\begin{array}{l}
d\omega_{\ }     =\omega \wedge \omega_{1}\\
d\omega_{1} =\omega \wedge \omega_{2}\\
\ \ \  \  \ \  \vdots\\
d\omega_{n} =\omega \wedge \omega_{n+1}+\sum_{k=1}^{n}\left(\begin{array}{c} n\\ k \end{array} \right)\omega_{k} \wedge \omega_{n-k+1}
\end{array}
$$

Elle sera dite de longueur $\ell$ si $\omega_{i}=0$ pour $i \geq \ell$, de longueur $1$ si il existe $F$ m\'{e}romorphe et un entier $k$ tels que 
$d(F^{1/k}\omega)=0$.
\end{definition}

On remarquera que les suites de longueur 1 sont des suites de longueur 2 particuli\`eres correspondant \`a $\omega_1=\frac{1}{k}\frac{dF}{F}$. La fonction 
multivalu\'ee $F^{1/k}$ est appel\'ee facteur int\'egrant de la forme $\omega$. Sur $\Delta$ une telle suite existe toujours gr\^ace \`a l'algorithme de 
Godbillon-Vey. Par contre l'existence d'une suite de longueur finie pour $\omega$ n'est pas toujours assur\'ee. La longueur minimale des suites, mais pas 
la suite elle-m\^{e}me, ne d\'{e}pend que du feuilletage $\mathcal{F}_\omega$. Nous ne nous int\'eresserons qu'aux suites de 
Godbillon-Vey de longueur inf\'erieur ou \'egale \`a trois. Pour un \'etude des suites de longueur sup\'erieure, on pourra consulter \cite{CLLPT}.
Les feuilletages de codimension un admettant une suite de Godbillon-Vey de 
longueur inf\'erieure ou \'egale \`a trois sont caract\'eris\'es par leur groupo\"ide de Galois :

\begin{theoreme}
\label{suitegv}
Le feuilletage $\mathcal{F}_{\omega}$ admet une suite de Godbillon-Vey de longueur $\ell$ avec $\ell \leq 3$ si et seulement si son groupo\"{i}de de Galois 
est contenu dans un $\mathcal{D}$-groupo\"ide de Lie de rang transverse $\ell$.
\end{theoreme}

\begin{preuved}{du th\'eor\`eme \ref{suitegv} pour $\ell =1$}
Si $\mathcal{F}$ admet une suite de Godbillon-Vey de longueur 1, il existe $F$ m\'eromorphe et un entier $k$ tels que $d(F^{1/k}\omega )=0$. Soit $\Gamma $
 un diff\'{e}omorphisme local conservant le feuilletage, c'est-\`{a}-dire v\'{e}rifiant $\Gamma ^{\ast }(\omega )=f_{\Gamma }\omega $ pour une fonction 
 $f_{\Gamma }$. On a alors $\Gamma^{\ast }(F^{1/k}\omega )=F^{1/k}\circ \Gamma \ f_{\Gamma }\omega $. Consid\'{e}rons l'\'{e}quation d'invariance de la 
 forme ferm\'ee :
$$
\Gamma^{\ast }(F^{1/k}\omega )=F^{1/k}\omega.
$$
En prenant la puissance $k$-i\`eme, nous obtenons l'\'equation \`a coefficients m\'eromorphes: 
$$
F=F\circ \Gamma \ f_{\Gamma }^k.
$$ 
Ceci est l'\'{e}quation d'un $\mathcal{D}$-groupo\"{i}de de Lie. En effet $f_{\Gamma}$ est un polyn\^ome en les d\'eriv\'ees premi\`eres de $\Gamma$ \`a
 coefficients holomorphes en $\Gamma$ et m\'eromorphes en $x$. De plus l'\'egalit\'e
$$
f_{\Gamma_1 \circ\Gamma_2}=(f_{\Gamma_{1}}\circ \Gamma_{2}) f_{\Gamma_{2}} 
$$
permet de v\'erifier les axiomes d'un $\mathcal{D}$-groupo\"{i}de de Lie :
\begin{multline}
\ \ \ F\circ (\Gamma_1 \circ \Gamma_2)\ f_{\Gamma_1 \circ  \Gamma_2}^k-F \nonumber \\
= \left( (F\circ \Gamma_1 \ f_{\Gamma_1 }^k-F)f_{\Gamma_2 }^k \right) \circ \Gamma_2-(F\circ \Gamma_2 \ f_{\Gamma_2 }^k-F). \ \ \ \nonumber
\end{multline}
V\'erifions que ce $\mathcal{D}$-groupo\"{i}de de Lie est admissible pour le feuilletage $\mathcal{F}_\omega$. Notons $L_X$ la d\'eriv\'e de Lie par rapport
 \`a un champ $X$ et pour un champ $X$ pr\'eservant le feuilletage d\'efinissons $f_{X}$ par $L_{X}\omega = f_{X}\omega$. La $\mathcal{D}$-alg\`ebre de Lie 
 du groupo\"ide d'invariance de la forme ferm\'ee a pour \'{e}quation:
$$
 L_{X}F+kFf_{X}=0 
$$
obtenue en lin\'earisant l'\'equation $F=F\circ \Gamma \ f_{\Gamma }^k$. Pour tout champ $X$ v\'{e}rifiant $\omega (X)=0$, on a
$$
\begin{array}{rl}

d(F^{1/k}\omega )(X,.)& =d(F^{1/k})(X)\omega + F^{1/k}d\omega(X,.)\\
                      & =(L_{X}F^{1/k}+F^{1/k}f_{X})\omega \\
                      & =(\frac{1}{k}\frac{L_{X}F}{F}+f_{X})F^{1/k}\omega.

\end{array} 
$$
L'\'equation $d(F^{1/k}\omega)=0$ est \'{e}quivalente au fait que tout champ de vecteur tangent au feuillatage v\'erifie l'\'equation de la 
$\mathcal{D}$-alg\`ebre de Lie. Par d\'efinition de $\mathcal{G}al(\mathcal{F}_{\omega})$, ce dernier est inclus dans le $\mathcal{D}$-groupo\"{i}de de Lie
 d'\'{e}quations :
$$
\Gamma^* \omega \wedge \omega = 0 \text{ et }F^k=F^k\circ \Gamma f_{\Gamma }^k.
$$
Cette \'equation \'etant d'ordre un, son expression locale donne une \'equation d'ordre un qui correspond donc \`a un \dgl\ de rang transverse un.

R\'{e}ciproquement, supposons que $\mathcal{F}$ admette un groupo\"{i}de de Galois de rang transverse un. Sur chaque carte de redressement de coordonn\'{e}e
 transverse $t$, son \'equation, donn\'ee par le (3) du th\'eor\`eme \ref{integr1} et le lemme \ref{local}, est celle du
 $\mathcal{D}$-groupo\"{i}de de Lie d'invariance de la forme m\'{e}romorphe $(\gamma (t)dt)^{\otimes k}$ o\`u $\gamma$ est une racine $k$-i\`eme de 
 $\gamma_k$. Par transitivité, cet entier $k$ est ind\'ependant de la carte locale. Les p\^oles et les z\'eros des diff\'erentes formes locales 
 $\gamma^{\otimes k}$ se recollent en un ensemble analytique $Z$ de codimension un dans $\Delta - Sing(\mathcal{F})$. Dans chaque carte de redressement du 
 feuilletage ne rencontrant pas $Z$, d\'efinissons le facteur int\'egrant local par 
$$ 
\Phi(z,t) = \frac{\gamma (t)}{w(z,t)} 
$$
o\`u  $w$ est d\'efini en coordonn\'ees par $\omega=  w(z,t)dt$.\\
Sur un autre ouvert muni de coordonn\'ees redressantes $(\widetilde{z},\widetilde{t})$, consid\'erons la forme $\widetilde{\gamma}$ et le facteur 
int\'egrant $\widetilde{\Phi}$ d\'efinis de mani\`ere analogue. 

Des changements de coordonn\'ees :
$$ 
\begin{array}{l}
\ \ \ \ \widetilde{\gamma}(\widetilde{t})\left(\frac{\partial \widetilde{t}}{\partial t}\right)=\zeta \gamma (t) \text{ avec } \zeta^{k}=1 \\ 
\text{ et } \widetilde{w}(\widetilde{z},\widetilde{t})\frac{\partial \widetilde{t}}{\partial t}=w(z,t),
\end{array}
$$
on d\'eduit que les facteurs int\'egrants $\Phi$ se recollent en dehors de $Z$ \`a multiplication par une racine $k$-i\`eme de l'unit\'e pr\`es.
La fonction $F=\Phi^k$ est bien d\'efinie et est m\'eromorphe en dehors du lieu singulier du feuilletage. Celui-ci est de codimension deux et 
le th\'eor\`eme de Hartogs assure son prolongement m\'eromorphe sur $\Delta$. Toute racine $k$-i\`eme de cette fonction est un facteur int\'egrant 
et d\'efinit une suite de Godbillon-Vey de longueur un pour le feuilletage.
  
\end{preuved}

\begin{remarque}
Si $\mathcal{F}$ n'a pas d'int\'egrale premi\`ere m\'eromorphe, la forme $F^{1/k} \omega$ est unique \`a multiplication par une constante pr\`es. 
Les \'{e}quations du $\mathcal{D}$-groupo\"{i}de que nous obtenons sont ind\'{e}pendantes de la forme $\omega$ initialement choisie. 
\end{remarque}

\begin{preuved}{du th\'eor\`eme \ref{suitegv} pour $\ell=2$}
Supposons que $\mathcal{F}$ admette une suite de Godbillon-Vey de longueur deux donn\'ee par une forme $\alpha$ telle que $d\omega = \omega \wedge \alpha $
 et $d\alpha=0$. Soit $(X,Y)=\Gamma(x,y)$ un germe de diff\'{e}omorphisme local pr\'{e}servant le feuilletage, c'est-\`a-dire v\'erifiant :
 
$$
\Gamma^{\ast}\omega =f_{\Gamma}\omega.
$$
Les \'{e}galit\'{e}s
$$
\Gamma^{\ast}d\omega = \Gamma^{\ast}\omega \wedge \Gamma^{\ast}\alpha = -f_{\Gamma}\Gamma^{\ast}\alpha \wedge \omega
$$
$$
d\Gamma^{\ast}\omega = df_{\Gamma} \wedge \omega + f_{\Gamma}d\omega = (df_{\Gamma}-f_{\Gamma} \alpha)\wedge \omega
$$
donnent la relation
$$
(\Gamma^{\ast}\alpha-\alpha + \frac{df_{\Gamma}}{f_{\Gamma}})\wedge\omega = 0
$$ 
entre les facteurs int\'egrants $\alpha$ pour $\omega$ et $\Gamma^* \alpha$ pour
$\Gamma^* \omega$.
Il existe donc une fonction $g_{\Gamma}$, d\'etermin\'ee par $\Gamma$, ses d\'eriv\'ees premi\`eres et secondes, v\'{e}rifiant 
$$
\Gamma^{\ast}\alpha-\alpha + \frac{df_{\Gamma}}{f_{\Gamma}} = g_{\Gamma}\omega. 
$$ 
On v\'erifie que les coefficients $g_{\Gamma}$ satisfont 
$$
g_{\Gamma_{1}\circ\Gamma_{2}} = (g_{\Gamma_{1}}\circ\Gamma_{2})f_{\Gamma_{2}} + g_{\Gamma_{2}}.
$$ 
Les transformations $\Gamma $ telles que $g_{\Gamma}=0 $ sont donc solutions d'un $\mathcal{D}$-groupo\"{i}de de Lie d'ordre deux. Montrons que ce 
$\mathcal{D}$-groupo\"{i}de est admissible. Plus pr\'ecisement, montrons qu'il contient $\mathcal{G}al(\mathcal{F}_{\omega})$ si et seulement si la forme 
$\alpha$ est ferm\'{e}e. En reprenant les notations du cas pr\'ec\'edent, les \'{e}quations lin\'{e}aris\'{e}es le long de l'identit\'{e} de $g_{\Gamma} = 0
$ sont 
$$ 
L_X\omega = f_X \omega  \text{ et } df_{X} + L_{X}\alpha = 0. 
$$ 
Prenons un champ $X$ tel que $\omega(X)=0$. On a alors:
$$
\begin{array}{lrl}

            & L_{X}\omega  & = d(\omega(X))+d\omega(X,.)\\
            &              & = -\alpha(X)\omega \\
\text{donc} & df_{X}       & = -d(\alpha(X)) \\
\text{et}   & L_{X}\alpha  & = d(\alpha(X))+d\alpha(X,.). 
\end{array}
$$
Nous en d\'eduisons que $X$ est solution du syst\`eme lin\'{e}aris\'{e} si et seulement si $d\alpha=0$. L'\'equation $g_{\Gamma}=0 $ \'etant d'ordre deux,
 le \dgl\ admissible que nous venons de construire est de rang transverse deux.\\
% REMardque plus haut rangtransverse=ordre de lequation jeute apres la def de
% rang transv
\indent R\'{e}ciproquement, supposons qu'il existe un \dgl\ admissible de rang transverse deux. Nous allons utiliser l'expression locale de ce \dgl\ pour 
construire une int\'egrale premi\`ere \`a monodromie affine. Pla\c{c}ons nous sur un ouvert de redressement de coordonn\'{e}es $(t,z)$ du feuilletage. 
La forme $\omega$ s'\'{e}crit $w(t,z)dt$ et les \'{e}quations du \dgl\ admissible, donn\'ees par lemme \ref{local} et le th\'eor\`eme 
\ref{integr1}, sont de la forme :
$$ 
\frac{\partial T}{\partial z}=0 \text{ et } \mu(t) = \mu(T)\frac{\partial T}{\partial t} + \frac{\frac{\partial^{2} T}{\partial t^2}}{\frac{\partial T}{\partial t}}. 
$$
Quitte \`a le restreindre, supposons que cet ouvert soit simplement connexe et ne contienne pas de p\^oles de $\mu$. Consid\'{e}rons sur cet ouvert une 
int\'egrale premi\`ere $H$ solution des \'equations :  
$$ 
\frac{\partial H}{\partial z}=0 \text{ et } \frac{\frac{\partial^{2}H}{\partial t^{2}}}{\frac{\partial H}{\partial t}}=\mu.
$$ 
Sur un ouvert analogue muni de coordonn\'{e}es $(\widetilde{t},\widetilde{z})$ nous construisons de m\^eme une int\'egrale premi\`ere $\widetilde{H}$. 
Un calcul direct de changement de coordonn\'ees dans les \'equations d'un \dgl\ sur un disque (\cite{C2}) donne 
$$\widetilde{\mu}(\widetilde{t})=\mu(t)\frac{\partial t}{\partial
 \widetilde{t}} + \frac{\frac{\partial^{2} t}{\partial
 \widetilde{t}^2}}{\frac{\partial t}{\partial \widetilde{t}}}.
$$
Ainsi $\widetilde{H}$ v\'erifie 
$$
\frac{\frac{\partial^{2} \widetilde{H}}{\partial
    \widetilde{t}^{2}}}{\frac{\partial  \widetilde{H}}{\partial
    \widetilde{t}}}= \frac{\frac{\partial^{2} \widetilde{H}}{\partial
    t^{2}}}{\frac{\partial  \widetilde{H}}{\partial t}}\frac{\partial
  t}{\partial \widetilde{t} } +  \frac{\frac{\partial^{2} t}{\partial
 \widetilde{t}^2}}{\frac{\partial t}{\partial \widetilde{t}}} = \widetilde{\mu},
$$ 
d'o\`u
$$
\frac{\frac{\partial^{2} \widetilde{H}}{\partial
    t^{2}}}{\frac{\partial  \widetilde{H}}{\partial
    t}}=\mu=\frac{\frac{\partial^{2} H}{\partial
    t^{2}}}{\frac{\partial  H}{\partial
    t}}.
$$
Sur l'intersection des deux ouverts, que l'on suppose connexe, il existe donc deux constantes $a$ et $b$ telles que $\widetilde{H}=aH+b$. En prolongeant 
une solution locale par cette formule, on construit une int\'{e}grale premi\`{e}re $H$ du feuilletage, en dehors du lieu des p\^oles des diff\'erents $\mu$, 
multivalu\'{e}e \`{a} monodromie affine. La fonction $F=\frac{\frac{\partial H}{\partial t}}{w}$ v\'{e}rifie 
$$
d(Fwdt)=0.% \text{ et } \frac{dF}{F} = \frac{\frac{\partial^{2} H}{\partial t^{2}}}{\frac{\partial H}{\partial t}}dt - \frac{dw}{w}.
$$
La forme ferm\'ee $\alpha= \frac{dF}{F}$ est univalu\'ee en dehors des p\^oles des diff\'erents $\mu$ et v\'erifie $d\omega = \omega \wedge \alpha$. Sur 
une carte de redressement contenant des p\^oles de $\mu$, la forme $\alpha = \frac{dF}{F}$ se prolonge m\'{e}romorphiquement. En effet, 
$$
\alpha =\frac{dF}{F} = \frac{\frac{\partial^{2} H}{\partial
    t^{2}}}{\frac{\partial H}{\partial t}}dt - \frac{dw}{w}=\mu(t)dt - \frac{dw}{w},
$$
o\`u les fonctions $\mu $ et $w$ sont m\'eromorphes sur chaque carte de redressement. Le lieu singulier du feuilletage \'etant de codimension au moins deux,
 la forme $\alpha$ se prolonge m\'eromorphiquement au polydisque.  
\end{preuved}

\begin{remarque}
Les suites de Godbillon-Vey de longueur deux de la forme $(\omega,\alpha)$ et $(f\omega, \alpha -\frac{df}{f})$ sont \'equivalentes  
(voir \emph{\cite{scardua}} et \emph{\cite{G}}) : elles d\'efinissent la m\^eme structure affine transverse en dehors de $Z$. S'il n'existe pas de 
facteur int\'egrant m\'eromorphe, la suite $(\omega,\alpha)$ est unique \`a \'equivalence pr\`es. L'\'equation du \dgl\ admissible de rang deux que
 nous obtenons est alors unique.
\end{remarque}

\begin{preuved}{du th\'eor\`eme \ref{suitegv} pour $\ell=3$} Nous suivrons la m\^{e}me strat\'{e}gie que dans les cas pr\'{e}c\'{e}dents.
 Supposons qu'il existe des formes m\'eromorphes $\alpha$ et $\beta$ telles que:
$$
\begin{array}{rl}
d\omega & =\omega \wedge \alpha \\
d\alpha & =\omega \wedge \beta \\
d\beta  & =\alpha \wedge \beta .
\end{array}
$$
Soit $\Gamma$ un automorphisme local du feuilletage et $f_\Gamma$ le coefficient de proportionnalit\'e qu'il d\'efinit. De la premi\`ere \'equation 
de Godbillon-Vey, nous d\'eduisons l'existence d'une fonction $g_{\Gamma}$, d\'etermin\'ee par $\Gamma$, ses d\'eriv\'ees premi\`eres et secondes, 
v\'{e}rifiant :
$$
\Gamma^{\ast}\alpha-\alpha + \frac{df_{\Gamma}}{f_{\Gamma}} = g_{\Gamma}\omega. 
$$
De la deuxi\`eme \'equation, nous d\'eduisons les \'egalit\'es:
$$
d(\Gamma^* \alpha)= \omega \wedge(f_\Gamma\ \Gamma^* \beta)
$$
$$
d(\alpha - \frac{df_\Gamma}{f_\Gamma}+g_\Gamma \omega)= \omega \wedge \beta +
dg_{\Gamma} \wedge \omega + g_{\Gamma}\ \omega \wedge \alpha.
$$
En faisant la diff\'erence, on obtient:
$$
(f_\Gamma\ \Gamma^* \beta - \beta + dg_{\Gamma}- g_{\Gamma}\alpha)\wedge \omega =0.
$$
Il existe donc une fonction $h_{\Gamma}$ d\'ependant des d\'eriv\'ees troisi\`emes de $\Gamma$ telle que 
$$
f_{\Gamma }\Gamma ^{\ast }\beta -\beta -g_{\Gamma }\alpha +dg_{\Gamma }=h_{\Gamma }\omega. 
$$
Nous en d\'eduisons : 
$$
h_{\Gamma _{1}\circ \Gamma _{2}}=h_{\Gamma _{2}}+h_{\Gamma _{1}}\circ \Gamma_{2}\ (f_{\Gamma _{2}})^{2}+g_{\Gamma _{1}}\circ \Gamma _{2}\ g_{\Gamma_{2}}\ f_{\Gamma _{2}}
$$ 
d'o\`u :
$$
h_{\Gamma _{1}\circ \Gamma _{2}}-\frac{1}{2}(g_{\Gamma _{1} \circ \Gamma_{2}})^{2}=(h_{\Gamma _{2}}-\frac{1}{2}(g_{\Gamma _{2}})^{2})+(h_{\Gamma_{1}}-\frac{1}{2}(g_{\Gamma _{1}})^{2})\circ\Gamma _{2}\ (f_{\Gamma _{2}})^{2}.
$$
L'\'{e}quation $h_{\Gamma }-\frac{1}{2}(g_{\Gamma })^{2}=0$ v\'erifie les axiomes d'un $\mathcal{D}$-groupo\"{i}de d'ordre trois. Montrons que ce 
$\mathcal{D}$-groupo\"{i}de contient $\mathcal{G}al(\mathcal{F}_{\omega})$ si et seulement si la troisi\`{e}me \'{e}quation de la suite de Godbillon-Vey 
est v\'{e}rifi\'{e}e. La $\D$-alg\`ebre de Lie de ce \dgl\ a pour \'equation: 
$$
f_{X}\beta + L_{X}\beta - g_{X}\alpha + dg_{X} = 0
$$ 
o\`u $g_{X}\omega = df_X + L_X \alpha$. Soit $X$ tel que $\omega(X) = 0$. Comme dans le preuve pr\'ec\'edente, on a $ f_{X}=-\alpha(X) $. Nous d\'eduisons 
$$ 
g_{X}\omega =-d(\alpha (X))+d(\alpha (X))+d\alpha (X,.)=-\beta (X)\omega. 
$$
Des \'egalit\'es
$$
d\beta (X,.)=L_X\beta +d(g_X)
$$
$$
\alpha \wedge \beta (X,.)=-f_X \beta +g_X \alpha 
$$
on obtient par diff\'erence l'\'{e}quation de la $\D$-alg\`ebre de Lie sous la forme : 
$$ 
(d\beta - \alpha \wedge \beta)(X,.)=0 .
$$ 
La troisi\`{e}me \'{e}quation de la suite de Godbillon-Vey est donc \'{e}quivalente au fait que tout champ tangent au feuilletage est solution de la
 $\mathcal{D}$-alg\`ebre de Lie du $\mathcal{D}$-groupo\"{i}de de Lie que nous venons de construire. Ceci prouve l'inclusion de 
 $\mathcal{G}al(\mathcal{F}_{\omega})$ dans un \dgl\ admissible de rang transverse trois.

R\'{e}ciproquement, supposons que le feuilletage admette un \dgl\ admissible de rang transverse trois. Sur un ouvert de redressement du feuilletage de
 coordonn\'{e}es $(t,z)$, les \'{e}quations du \dgl\ admissible, donn\'ees par lemme \ref{local} et le th\'eor\`eme \ref{integr1}, 
 sont de la forme :  
$$ 
\frac{\partial T}{\partial z}=0 \text{ et } \nu(t) = \nu(T)\left(\frac{\partial T}{\partial t}\right)^{2}+ S_{t}(T)
$$
o\`{u} $ S_{t}(T)= 2\frac{\frac{\partial^3 T}{\partial t^3}}{\frac{\partial T}{\partial t}} -3 \left(\frac{\frac{\partial^2 T}{\partial t^2}}{\frac{\partial T}{\partial t}}\right)^2$ 
est la schwartzienne de $T$ par rapport \`{a} $t$. Nous allons nous servir de $\nu$ pour construire une int\'{e}grale premi\`{e}re du feuilletage $H$ 
\`a monodromie projective et un couple de formes m\'eromorphes $(\alpha, \beta)$ v\'erifiant les \'equations de Godbillon-Vey. Nous pouvons toujours 
choisir une forme m\'{e}romorphe $\alpha$ v\'erifiant la premi\`{e}re \'{e}quation : il suffit de prendre un champ m\'{e}romorphe $X$ v\'{e}rifiant 
$\omega(X)=1$ et de poser $\alpha=L_{X}\omega$. Pla\c{c}ons-nous sur un ouvert de redressement ne rencontrant pas le lieu des p\^oles $Z$ des diff\'erents 
$\nu$. Soit $H$ une int\'egrale premi\`ere sur cet ouvert solution des \'equations: 
$$
\frac{\partial H}{\partial z}=0 \text{ et } S_{t}H=\nu(t). 
$$
Soit $\widetilde{H}$ une int\'egrale premi\`ere construite de mani\`ere analogue sur un ouvert de redressement muni des coordonn\'ees 
$(\widetilde{t},\widetilde{z})$. D'apr\`es les changements de variables usuels sur les d\'eriv\'ees Schwartziennes et sur les \'equations 
des \dgl\ \cite{C2}, nous avons
$$
S_t \widetilde{H}=S_{\widetilde{t}} \widetilde{H} \left(\frac{\partial t}{\partial
    \widetilde{t}}\right)^2 + S_{\widetilde{t}} t
$$
et
$$
\nu(t)=\widetilde{\nu}({\widetilde{t}}) \left(\frac{\partial t}{\partial
    \widetilde{t}}\right)^2 + S_{\widetilde{t}} t.
$$
Nous en d\'eduisons que $S_t \widetilde{H}=S_t H$ et donc que $H$ se prolonge de mani\`{e}re multivalu\'{e}e sur le compl\'ementaire de $Z$ avec
 une monodromie projective. A partir de cette int\'egrale premi\`ere, on construit la fonction $F$ par $ dH=F\omega $ d'o\`u :
$$
\frac{dF}{F} \wedge \omega + d\omega = 0.
$$
Contrairement au cas pr\'ec\'edent, la forme $\frac{dF}{F}$ n'est pas m\'eromorphe. La forme $\gamma = \frac{dF}{F}- \alpha $ v\'{e}rifie 
$\gamma = G\omega $ pour une certaine fonction $G$. On a alors $d\alpha = \omega \wedge (dG-G\alpha)$. La forme $\beta$ cherch\'ee est de la forme : 
$$ 
\beta = dG - G\alpha + K\omega .
$$ 
En rempla\c{c}ant cette expression dans la troisi\`{e}me \'equation de Godbillon-Vey, on obtient :
$$
(dK+GdG) \wedge \omega + (G^2 +2K)\omega \wedge \alpha =0.
$$
Posons $K=-\frac{1}{2}G^{2}$.
Ceci nous permet de construire la forme $\beta$ \`a partir de $\alpha$ et $G$ v\'erifiant les \'equations de Godbillon-Vey. Il nous reste \`{a} montrer
 sa m\'{e}romorphie. En prenant une autre d\'etermination $\widetilde{H}=\frac{aH+b}{cH+e}$, nous obtenons, \textit{a priori}, une autre forme 
 $\widetilde{\beta}$. En calculant cette forme, on a : 

$$
\begin{array}{rcl}

d\widetilde{H}=\widetilde{F}\omega   &  =   &\frac{F}{(cH+e)^{2}}\omega \\
\frac{d\widetilde{F}}{\widetilde{F}} &  =   &\frac{dF}{F}-2\frac{c}{cH+e}dH \\
\widetilde{G}                       &  =   &G-\frac{2cF}{cH+e}\\
d\widetilde{G}                      &  =   &dG-\frac{2cdF}{cH+e}+\frac{2c^{2}F^{2}\omega}{(cH+e)^{2}}\\
\widetilde{G}\alpha                 &  =   &G\alpha-\frac{2cF\alpha}{cH+e}\\
\frac{1}{2}\widetilde{G}^{2}\omega  &  =   &\frac{1}{2}G^{2}\omega -\frac{2cGF \omega}{cH+e}+\frac{2c^{2}F^{2}\omega}{(cH+e)^{2}}.

\end{array}
$$
En rempla\c{c}ant $dF$ par $F \alpha +FG \omega $, on v\'erifie que $\beta=\widetilde{\beta}$. La forme $\beta$ est donc univalu\'ee sur les ouverts 
de redressement ne rencontrant pas $Z$. V\'erifions qu'elle admet un prolongement m\'eromorphe sur $Z$. Sur une carte de redressement on \'ecrit 
$\omega = w(t,z)dt$ et $\alpha= a_t(t,z)dt + a_z(t,z)dz$. Dans ces coordonn\'ees $F=\frac{1}{w}\frac{\partial H}{\partial t}$ d'o\`u :
$$
\begin{array}{rcl}
Gw dt & = & \left(\frac{\partial_{t,t} H}{\partial_t H}-\frac{\partial_t
    w}{w}-a_t\right)dt + \left(\frac{\partial_z w}{w}+a_z \right)dz \\
G & = & \frac{1}{w}\left(\frac{\partial_{t,t} H}{\partial_t H}-\frac{\partial_t w}{w}-a_t\right)  \text{ et } \frac{\partial_z w}{w}+a_z=0 \\
dG & = & -\frac{\partial_t w}{w}Gdt-\frac{\partial_z
    w}{w}Gdz+\frac{1}{w}\left[\partial_t\left(\frac{\partial_{t,t}
    H}{\partial_t H}-\frac{\partial_t w}{w}-a_t\right)\right]dt \\
 & & +\frac{1}{w}\left(\partial_z\left(-\frac{\partial_t w}{w}-a_t\right)\right)dz \\
-G\alpha & = & -a_tGdt-a_zGdz \\
-\frac{1}{2}G^2\omega & = & \frac{1}{w}\left[-\frac{1}{2} \left(\frac{\partial_{t,t} H}{\partial_t H}\right)^2 - \frac{1}{2} \left(\frac{\partial_t w}{w} \right)^2 \right.\\ 
&& \hfill \left. -\frac{1}{2}a_t^2 + \frac{\partial_{t,t} H}{\partial_t H}\frac{\partial_t w}{w} + \frac{\partial_{t,t} H}{\partial_t H}a_t - \frac{\partial_t w}{w}a_t  \right].
 
\end{array}
$$
En sommant les trois derni\`eres \'equations, apr\`es simplification, on trouve :
$$
\beta = \frac{1}{w}\left[\nu(t)dt-d\left(\frac{\partial_t w}{w}+a_t\right)+\frac{1}{2}\left(\frac{\partial_t w}{w}+a_t\right)^2dt\right].
$$
On en d\'eduit que la forme $\beta $ est m\'{e}romorphe en dehors du lieu singulier du feuilletage et se prolonge m\'{e}romorphiquement \`a celui-ci.
\end{preuved}

\begin{remarque}
Les suites de Godbillon-Vey de longueur trois de la forme $(\omega,\alpha,\beta)$ et $(f\omega, \alpha -\frac{df}{f} + g\omega , \frac{1}{f}(\beta -dg +g\alpha + \frac{g^2}{2}\omega))$ 
sont \'equivalentes (voir \emph{\cite{scardua}} et \emph{\cite{G}}) : elles d\'efinissent la m\^eme structure transverse projective en dehors de $Z$. 
Dans le cas o\`u le feuilletage n'admet pas de suite de Godbillon-Vey de longueur deux, la suite de longueur trois est unique \`a \'equivalence pr\`es. 
Le \dgl\ obtenu est ind\'ependant de la suite.
\end{remarque}

\begin{remarque}
\label{formules}
Au cours de la preuve de ce th\'eor\`eme, nous avons donn\'e les  \'equations explicites d'un \dgl\ admissible pour le feuilletage ainsi que son 
expression locale sur une transverse. Soit $(t,z)$ des coordonn\'ees de redressement. La forme $\omega$ s'\'ecrit $w(t,z)dt$. L'expression locale 
du \dgl\ transverse est donn\'ee par les formules suivantes : 
\begin{enumerate}
\item Soit $(\omega, F^{1/k})$ une suite de longueur un pour le feuilletage. Le \dgl\ transverse est $\G_1(\mu)$ avec :
$$
\mu(t)dt=\frac{1}{k}\frac{dF}{F}+\frac{dw}{w}.
$$  
\item  Soit $(\omega, \alpha)$ une suite de longueur deux pour le feuilletage. Le \dgl\ transverse est $\G_2(\mu)$ avec :
$$
\mu(t)dt=\alpha+\frac{dw}{w}.
$$ 
\item  Soit $(\omega, \alpha, \beta)$ une suite de longueur trois pour le feuilletage. La forme $\alpha$ s'\'ecrit $a_tdt+a_zdz$. Le \dgl\ transverse est $\G_3(\nu)$ avec :
$$
\nu(t)dt=w \beta +d\left(  \frac{\partial_t w}{w}+a_t\right)- \frac{1}{2}\left(  \frac{\partial_t w}{w}+a_t\right)dt.
$$
R\'eciproquement, ces m\^emes formules permettent d'obtenir une suite de Godbillon-Vey explicite \`a partir des \'equations d'un \dgl\ admissible. 
\end{enumerate}   
\end{remarque}

\section{Groupo\"ides de Galois et int\'egrales premi\`eres }
\label{transcendance}

Rappelons les types de transcendances d'extensions du corps des fonctions m\'eromorphes sur un polydisque $\Delta$ de $\CC^n$.

\begin{definitions}\text{  \\ }
\label{deftypedetransc}
\begin{enumerate}
\item Une extension diff\'{e}rentielle du corps des fonctions m\'eromorphes
  sera dite de type Darboux si elle est obtenue par une suite d'extensions qui
  sont soit alg\'ebriques soit du type $ K(G) \supset K$ avec $dG=\gamma$, $\gamma$ \'etant une forme \`a coefficients dans $K$.
\item Une extension diff\'{e}rentielle du corps des fonctions m\'eromorphes sera dite Liouvillienne si elle est obtenue par une suite d'extensions qui
  sont soit alg\'ebriques soit du type $K(G) \supset K$ avec $dG=G\gamma_{1}+\gamma_{0}$, $\gamma_{1}$ et $\gamma_{0}$ \'etant des formes \`a coefficients 
  dans $K$.
\item Une extension diff\'{e}rentielle du corps des fonctions m\'eromorphes sera dite de type Riccati si elle est obtenue par une suite d'extensions qui
  sont soit alg\'ebriques soit du type $ K(G) \supset K$ avec $dG=\frac{G^{2}}{2}\gamma_{2}+G\gamma_{1}+\gamma_{0}$, $\gamma_{2}$, $\gamma_{1}$ et 
  $\gamma_{0}$ \'etant des formes \`a coefficients dans $K$.
\end{enumerate}
\end{definitions}

\begin{theoreme}
\label{typedetransc}
Soit $\mathcal{F}$ un germe de feuilletage de codimension un de $(\mathbb{C}^n,0)$. Le feuilletage $\mathcal{F}$ admet un int\'egrale premi\`ere 
m\'eromorphe (resp. de type Darboux, Liouville ou Riccati) si et seulement si $\mathcal{F}$ admet un \dgl\ admissible de rang transverse 0 (resp. 1,2 ou 3).
\end{theoreme}

Nous allons commencer par prouver le cas non transitif (rang transverse 0). Les trois autres affirmations seront prouv\'ees simultan\'ement par la suite.

\begin{lemme}
\label{principal}
Soient $\mathcal{F}$ un feuilletage de $\Delta$ de codimension un et $\mathcal{I}$ l'id\'eal du groupo\"ide de Galois de $\mathcal{F}$. 
Si le groupo\"ide de Galois n'est pas transitif alors l'id\'eal des \'equations d'ordre z\'ero, $\mathcal{I}_0=\mathcal{I} \cap \O_{J^*_0(\Delta)}$ 
de l'anneau $\O_{J^*_0(\Delta)}$, est engendr\'e par une unique \'equation. 
\end{lemme}

\begin{preuve}
Pla\c{c}ons-nous au voisinage d'un point r\'egulier du feuilletage. D'apr\`es la forme des \'equations locales du \dgl\ (voir le lemme \ref{local}), 
et le fait que sa $\D$-alg\`ebre de Lie soit de rang transverse nul, l'id\'eal $\I$ est engendr\'e par une \'equation d'ordre 0 (voir le (4) du 
th\'eor\`eme \ref{integr1}). Ceci signifie que l'id\'eal $\mathcal{I}$ est engendr\'e au voisinage de tout point de l'identit\'e 
$\{(x,x,id,0 \ldots ,0) | x \notin Sing(\mathcal{F})\}$ dans $J^*(\Delta)$ par une \'equation d'ordre z\'ero. Nous allons \'etendre cette propri\'et\'e 
\`a tout jet dont la source et le but en dehors d'un ensemble de codimension un.\\
En utilisant le th\'eor\`eme d'involutivit\'e g\'en\'erique pour les $\mathcal{D}$-groupo\"ides de Lie (th\'eor\`eme \ref{theodinvol}) et le 
th\'eor\`eme de Cartan-K\"ahler (\ref{theodeCK}), il existe un sous-ensemble analytique $Z$ de $\Delta$ et un entier $\ell$ tels que par tout point 
$a=(s(a), t(a),\ldots)$ de $J^*_\ell(\Delta)$ solution de $\mathcal{I}_\ell$ de source et but hors de $Z$ passe une solution convergente $\varphi$ de 
$\mathcal{I}_\ell$. Quitte \`a supposer $\ell$ assez grand, $\varphi$ est solution de $\mathcal{I}$ : en effet d'apr\`es le th\'eorème \ref{malg-ritt-rad},
 il existe un entier $\ell$ tel que $\mathcal{I}_\ell$ engendre diff\'erentiablement $\mathcal{I}$.

Par composition \`a la source, cette solution donne un isomorphisme d'un voisinage de $(t(a),t(a),id,0,\ldots)$ sur un voisinage de $a$ dans l'espace 
des jets d'ordre $\ell$. Puisque $\varphi$ est solution de $\I$ et que les z\'eros de $\I_\ell$ sont stables par composition en dehors de $Z$, cet 
isomorphisme se restreint en un isomorphisme local des espaces d\'efinis par $\mathcal{I}_\ell$ aux voisinages de ces m\^emes points. Nous en d\'eduisons 
qu'au voisinage de n'importe quel point au-dessus de source et but en dehors de $Z$, $\mathcal{I}$ est engendr\'e par une \'equation d'ordre 0. Le lieu des
 z\'eros $V$ de $\mathcal{I}_0$ est de codimension un dans $(\Delta-Z) \times (\Delta-Z)$. De plus cet ensemble analytique $V$ n'a pas de composante 
 irr\'eductible incluse dans $(Z \times \Delta) \cup (\Delta \times Z)$. Dans le cas contraire il existerait une fonction $f$ holomorphe sur $V$ nulle sur
 le compl\'ementaire de la composante irr\'eductible et non nulle sur celle-ci. Puisque le produit de cette fonction par une \'equation de $Z$ est 
 identiquement nul sur $V$, $f$ serait une \'el\'ement de torsion de $\O_{J^*(\Delta)}/\I$ pour une des deux projections. L'id\'eal $\I$ \'etant diff\'erentiel et r\'eduit, on v\'erifie que $f$ doit \^etre identiquement nulle. Le lieu des z\'eros de l'id\'eal r\'eduit $\mathcal{I}_0$ \'etant de codimension un, il est donc engendr\'e 
 par une \'equation $H$ (\cite{eisenbud}).  
\end{preuve}

Une relation d'\'equivalence analytique sur $\Delta$ est la donn\'ee d'un id\'eal $I$ de $\O_{\Delta \times \Delta}$ qui s'annule sur la diagonale, 
qui est stable par la sym\'etrie par rapport \`a la diagonale et qui v\'erifie la relation de transitivit\'e suivante :  
$$pr_{2,3}^*\I \subset pr_{1,2}^*\I + pr_{1,3}^*\I$$
 o\`u les $pr_{i,j}$ d\'esignent les trois projections naturelles de $\Delta \times \Delta \times \Delta$ sur $ \Delta \times \Delta $.

\begin{lemme}Sous les hypoth\`eses du lemme \ref{principal}, il existe un sous-ensemble analytique $Z$ de $\Delta$ tel que l'id\'eal 
$\mathcal{I}_0=(H(x,y))$ d\'efinisse une relation d'\'equivalence analytique sur $\Delta - Z$.
\end{lemme}

\begin{preuve}
L'id\'eal $\mathcal{I}_0$ \'etant form\'e des \'equations d'ordre z\'ero de l'id\'eal $\mathcal{I}$ d\'ecrivant un \dgl, les propri\'et\'es de 
r\'eflexivit\'e et de sym\'etrie sont v\'erifi\'ees. La stabilit\'e par composition nous donne l'inclusion de $pr_{2,3}^*\I_0$ dans l'id\'eal 
diff\'erentiablement engendr\'e par $pr_{1,2}^*\I_0 + pr_{1,3}^*\I_0$. Il nous faut v\'erifier qu'il est inclus dans l'id\'eal alg\'ebriquement 
engendr\'e par $pr_{1,2}^*\I_0 + pr_{1,3}^*\I_0$. Soit $Z$, l'ensemble analytique en dehors duquel on a la stabilit\'e par composition du \dgl\ 
(voir (3) de la d\'efinition \ref{defDGL}). Pla\c{c}ons-nous sur $(\Delta-Z) \times (\Delta-Z) \times (\Delta-Z)$ et consid\'erons les \'equations 
$pr_{1,2}^*H = H(x,y)$ et $pr_{1,3}^*H=H(x,z)$. Quitte \`a augmenter $Z$, les formes verticales pour la premi\`ere projection :
$$
\sum \frac{\partial H}{\partial y_i}(x,y) dy_i\  \text{ et }\  \sum \frac{\partial H}{\partial z_i}(x,z) dz_i
$$   
ne s'annulent pas.  On note $\mathcal{I}$ l'id\'eal diff\'erentiel engendr\'e par $(H(x,y),H(x,z))$ dans $\mathcal{O}_{J(\Delta \to \Delta \times \Delta)}$ 
et $\mathcal{I}_k$ les \'equations de $\mathcal{I}$ d'ordre inf\'erieur ou \'egal \`a $k$. Le fait que ces formes soient non nulles et non colin\'eaires 
nous permet d'utiliser une g\'enéralisation du th\'eor\`eme des fonctions implicites (voir \cite{tougeron}) :  pour tout z\'ero $(x,y,z)$ de $(H(x,y),H(x,z))$
 dans $(\Delta-Z) \times (\Delta-Z) \times (\Delta-Z)$, on peut trouver un z\'ero de $\I$ au-dessus de celui-ci pour la projection 
 $J_k(\Delta \to \Delta \times \Delta) \to \Delta \times \Delta \times \Delta$. Ceci signifie que $\mathcal{I}_{k} \cap \mathcal{O}_{J_0(\Delta \to \Delta \times \Delta)}$
co\"incide avec l'id\'eal alg\'ebriquement engendr\'e par $H(x,y)$ et $H(x,z)$ sur $(\Delta-Z) \times (\Delta-Z) \times (\Delta-Z)$ pour tout entier $k$. 
Les \'equations d'ordre z\'ero appartenant \`a $\mathcal{I}$ sont donc exactement celles de cet id\'eal.\\
Pour $\ell$ assez grand, $\mathcal{I}_\ell$ est un sous-groupo\"ide de $J^*_\ell(\Delta-Z)$, donc $H(y,z)$ appartient \`a $\mathcal{I}_\ell$. 
D'apr\`es ce qui pr\'ec\`ede  $H(y,z)$ appartient \`a l'id\'eal engendr\'e alg\'ebriquement par $H(x,y)$ et $H(x,z)$. Nous avons donc une relation 
d'\'equivalence en dehors de $Z$.
 \end{preuve}

%% \begin{remarque}
%% L'abscence de torsion implique que l'ensemble des $x$ tels que $H(x,y)$ soit identiquement nul est de codimension au moins deux. Mais il se peut qu'il soit strictement inclus dans le lieu singulier du feuilletage
%% \end{remarque}

\begin{preuved}{du cas m\'eromorphe du th\'eor\`eme \ref{typedetransc}}
On suppose que le groupo\"ide de Galois est d'ordre 0. Soit $H$ une \'equation
de l'id\'eal $\I_0$ donn\'ee par le lemme \ref{principal} et $Z$ le sous
ensemble analytique en dehors duquel on a la stabilit\'e du \dgl\ par
composition. Quitte \`a agrandir $Z$, nous supposerons qu'il contient le lieu singulier du feuilletage.
On note $R$ le lieu des z\'eros de $H$ dans $(\Delta \times \Delta)$ et $R|_{\Delta-Z}$ sa restriction sur $(\Delta-Z) \times( \Delta-Z)$.
Montrons que le quotient de $(\Delta-Z)$ par $R|_{\Delta-Z}$ est une surface de Riemann.\\
Par transitivit\'e, les classes d'\'equivalence $pr_2(pr_1^{-1}(p) \cap
R|_{\Delta-Z})$ sont constantes le long des feuilles du feuilletage. Les
projections $R|_{\Delta-Z}$ sur $\Delta-Z$ \'etant sans torsion, elles sont plates au-dessus d'une transverse en $p$ au
feuilletage. Par transitivit\'e elles sont plates sur un ouvert contenant
$p$. D'apr\`es \cite{F} elles sont ouvertes et en particulier le satur\'e pour
$R|_{\Delta-Z}$ d'un ouvert est un ouvert.

Pour prouver la s\'eparabilit\'e du quotient, on prend deux points $p$ et $q$
non \'equivalents. Soit $T$ une transverse au feuilletage en $p$. Le point (4) du
th\'eor\`eme \ref{integr1} nous assure que les classes d'\'equivalence de $p$
et de $q$ intersectent $T$ en un nombre fini de points. On peut donc s\'eparer
ces deux ensembles par des ouverts dans $T$ satur\'es pour la relation
d'\'equivalence restreinte \`a $T$. Les satur\'es de ces ouverts donnent deux
ouverts dans le quotient qui ne s'intersectent pas. Le quotient par
$R|_{\Delta-Z}$ est un espace topologique s\'epar\'e.

La construction d'un atlas holomorphe de cartes locales sur cet espace se fait de la
mani\`ere suivante. Au voisinage $U$ d'un point $p$ de $\Delta-Z$, il existe
une fonction holomorphe $h$ constante sur les classes d'\'equivalence. On
prolonge $h$ sur le satur\'e $\overline{U}=pr_2(pr_1^{-1}(U) \cap
R|_{\Delta-Z})$ du voisinage par $pr_{2_*}(pr_1^*(h)|_{ R|_{\Delta-Z}})$. Ceci
nous d\'efinit une carte sur l'ouvert $\overline{U}$ du quotient. Soient
$(h_1,\overline{U_1})$ et $(h_2,\overline{U_2})$ deux cartes d'intersection
non vide. Les applications $h_1|_{\overline{U_1}\cap \overline{U_2}}$ et
$h_2|_{\overline{U_1}\cap \overline{U_2}}$ ont les m\^emes hypersurfaces de
niveau. Il existe une application holomorphe $F$ telle que
$h_1|_{\overline{U_1}\cap \overline{U_2}}=F\circ h_2|_{\overline{U_1}\cap
  \overline{U_2}}$. Celle-ci d\'efinit un changement de carte holomorphe pour
la vari\'et\'e quotient. \\
On note alors $S$ la surface de Riemann obtenue et $\pi : \Delta-Z \to S$ le passage au quotient.
Montrons que quitte \`a rajouter des points \`a $S$, $\pi$ se prolonge \`a
$\Delta-Sing{\mathcal{F}}$. Si une composante irr\'eductible de $Z$ est
transverse au feuilletage, $\pi$ \'etant constante sur les feuilles, elle se
prolonge \`a cette composante. Sinon, consid\'erons une transverse $T$ au
feuilletage en un point $p$ de cette composante. D'apr\`es le (4) du th\'eor\`eme \ref{integr1}, il existe une
coordonn\'ee source $x$ sur $T$ et but $X$ sur $T$ dans laquelle
l'\'equation $H$ sur $T \times T$ s'\'ecrive $X^k-x^k=0$. Par construction de $S$, au voisinage de $p$ le passage au
quotient est donn\'e par :
$$\begin{array}{cc}
\pi: & T-\{p\} \to S \\
 & \ \ \ \ \ \ \ \ \ \ x \mapsto x^k.
\end{array}$$
Le quotient $\pi$ admet donc un prolongement holomorphe sur $Z$.\\
Notons encore $S$ l'image de ce prolongement. En ramenant une fonction m\'eromorphe de $S$ sur $\Delta - Sing \F$, 
on obtient une int\'egrale premi\`ere m\'eromorphe qui se prolonge \`a $\Delta$. 
R\'eciproquement si $\mathcal{F}$ admet une int\'egrale premi\`ere
m\'eromorphe, le groupo\"ide d'invariance de celle-ci est un \dgl\ admissible
pour $\mathcal{F}$ d'ordre z\'ero : voir l'exemple \ref{mero}.         
\end{preuved}

Pr\'ecisons maintenant la condition n\'ecessaire et suffisante sur
l'\'equation d'ordre z\'ero $H$ du groupo\"ide de Galois pour que le
feuilletage admette une int\'egrale premi\`ere holomorphe.   
\begin{proposition}
Lorsque le groupo\"ide de Galois de $\mathcal{F}$ est non transitif d'\'equation $H$,
le feuilletage admet une int\'egrale premi\`ere holomorphe si et seulement si $H(0,y)$ est non identiquement nulle.
\end{proposition}

\begin{preuve}
Si $H(0,y)$ est non identiquement nulle, nous pouvons supposer que $H$ est non identiquement nulle le long de l'axe des $y_n$ et appliquer le th\'eor\`eme de pr\'eparation de Weierstrass afin d'\'ecrire :
$$
H(x,y)=y_n^k+a_{k-1}(x,\overline{y})y_n^{k-1}+ \ldots +a_0(x,\overline{y})
$$
o\`u $\overline{y}=(y_1,\ldots,y_{n-1})$. Fixons deux points $x$ et $y$ en
dehors de $Z$ tels que $H(x,y)=0$.  
En utilisant la transitivit\'e de la relation d'\'equivalence, $H(x,y)=0$
implique qu'au voisinage de tout $z \in \Delta-Z$ il existe une unit\'e $u(z)$
telle que $H(x,z)=u(z)H(y,z)$. Gr\^ace aux normalisations de Weierstrass des
polyn\^omes $H(x,z)$ et $H(y,z)$ on obtient
$a_0(x,\overline{z})=a_0(y,\overline{z})$ pour tout $\overline{z}$. En
particulier on a $a_0(x,0)=a_0(y,0)$. Cette fonction est non constante. En
effet par sym\'etrie, il existe une unit\'e $v$ telle que $H(x,y)=v(x,y)H(y,x)$
d'o\`u 
$$
a_0(x,0)=H(x,0)=v(x,0)H(0,x)=v(x,0)(x_n^k + \ldots).
$$
Les feuilles de $\mathcal{F}$ \'etant incluses dans les classes
d'\'equivalence de la relation d'\'equivalence donn\'ee par $H$, la fonction
holomorphe $a_0(x,0)$ est une int\'egrale premi\`ere du feuilletage.\\
R\'eciproquement, si le feuilletage admet une int\'egrale premi\`ere
holomorphe non constante $h(x)$, l'\'equation $h(x)-h(y)$ definit un
$\mathcal{D}$-groupo\"ide de Lie contenant le groupo\"ide de Galois du
feuilletage. Les z\'eros de $H(x,y)$ sont donc inclus dans ceux de
$h(x)-h(y)$. En particulier pour $x=0$ ceci montre que $H(0,y)$ est non
identiquement nulle. 
\end{preuve}

\begin{remarque}
Cette preuve s'adapte au cas m\'eromorphe, suivant les indications de B. Malgrange, en utilisant le th\'eor\`eme de Weiertrass o\`u on consid\`ere les variables $x$ comme param\`etres.
\end{remarque}

Les preuves des autres cas du th\'eor\`eme \ref{typedetransc} (rang transverse
1, 2 ou 3) se d\'eduisent du th\'eor\`eme \ref{suitegv} et du th\'eor\`eme
ci-dessous :
\begin{theoreme}
\label{theosinger}
Un germe de feuilletage holomorphe singulier $\mathcal{F}$ de codimension
un admet une int\'egrale premi\`ere de type Darboux (resp. Liouville ou Riccati)
si et seulement s'il admet une suite de Godbillon-Vey de longueur un
(resp. deux ou trois).
\end{theoreme}
\begin{preuve}
Le cas Liouvillien est d\^u \`a M. Singer : \cite{S}. Sa g\'en\'eralisation au
cas Riccati est faite dans \cite{C}. Nous allons donner la preuve du cas Darboux.

Supposons qu'il existe une int\'egrale premi\`ere $H$ de type Darboux pour la forme $\omega$. On note
$K \subset K_1 \ldots \subset K_n$ la suite des extensions du corps $K$, telle que $H$ soit dans $K_n$. On supposera que cette suite est de longueur minimale parmi toutes les suites d'extensions de type Darboux n\'ecessaires \`a la construction d'une int\'egrale premi\`ere.\\
\indent L'extension $K_{n-1} \subset K_{n}$ ne peut pas \^etre alg\'ebrique. Dans le cas contraire, $H$ serait alg\'ebrique sur $K_{n-1}$. On note
$P(X)=X^p+a_{p-1}X^{p-1}\ldots+a_0$ son polyn\^ome minimal ; $H$ n'\'etant pas
constante, il en serait de m\^eme pour au moins un des $a_i$. On aurait
$$
0=dP(H)\wedge \omega = \sum_i H^i da_i \wedge \omega
$$
d'o\`u par minimalit\'e de $P$, $da_i\wedge\omega=0$. L'existence des
int\'egrales premi\`eres $a_i$ dans $K_{n-1}$ contredit la minimalit\'e de la suite d'extensions. La derni\`ere extension est donc transcendante. \\
\indent Soit $G$ telle que $K_n = K_{n-1}(G)$
o\`u $dG=\gamma$ est une forme \`a coefficients dans $K_{n-1}$. Lorsqu'on
\'ecrit $dH=F\omega$ dans $K_{n}$, on peut supposer que le facteur int\'egrant $F$ est dans
$K_{n-1}$. En effet, en \'ecrivant $F=a_kG^k+\ldots$ par division suivant les
puissances croissantes de $G$ et en calculant $d(F\omega)$,  par transcendance
de $G$ on obtient $d(a_k\omega)=0$. On peut donc consid\'erer la suite de longueur minimale donn\'ee par les $K_i$ pour $i$ inf\'erieur \`a $n-1$ et $K_n=K_{n-1}(H)$ avec $dH=a_k\omega$.
Il existe donc un facteur int\'egrant pour $\omega$ dans l'avant-dernier corps
de la suite d'extension donnant une int\'egrale premi\`ere de type Darboux.\\
\indent Si l'extension $K_{n-1}$ de $K_{n-2}$ est purement transcendante, le
m\^eme raisonnement assure l'existence d'un facteur int\'egrant pour $\omega$
dans $K_{n-2}$, ce qui contredit la minimalit\'e de la suite. L'extension
$K_{n-1}$ de $K_{n-2}$ est donc alg\'ebrique.\\
\indent Soit $F$ un facteur int\'egrant de $\omega$ dans $K_{n-1}$. Il est
alg\'ebrique sur $K_{n-2}$. On note $F_1,\ldots, F_p$ ses quantit\'es
conjugu\'ees. Comme $d(F_i\omega)=0$, on a 
$$
\frac{d(F_1\ldots F_p)}{F_1\ldots F_p}\wedge\omega=pd\omega.
$$
Le produit $\widetilde{F}= F_1\ldots F_p $ est un \'el\'ement non nul de
$K_{n-2}$ dont une racine $p$-i\`eme $\sqrt[p]{\widetilde{F}}$ est un facteur
int\'egrant pour $\omega$. On construit donc une nouvelle suite de longueur
minimale en conservant les $n-2$ premi\`eres
extensions et en rempla\c{c}ant $K_{n-1}$ par
$K_{n-2}(\sqrt[p]{\widetilde{F}})$ et $K_n$ par $K_{n-2}(\sqrt[p]{\widetilde{F}},H)$ avec
$dH=\sqrt[p]{\widetilde{F}}\omega$.
Montrons que cette suite est de longueur deux, c'est-\`a-dire $K_{n-2}=K_0$.\\
\indent Si l'extension $K_{n-2}$ de $K_{n-3}$ est alg\'ebrique, le
raisonnement pr\'ec\'edent permet de construire une fonction $\widetilde{F}$
dans $ K_{n-3}$ dont une racine est un facteur int\'egrant pour $\omega$. On
pourrait alors construire une suite de longueur $n-1$ qui contredirait la minimalit\'e de la suite.\\
\indent Si l'extension $K_{n-2} = K_{n-3}(G)$ est transcendante avec $dG$ \`a
coefficients dans $K_{n-3}$. On \'ecrit
$\widetilde{F}= a_kG^k(1+ a_1 G^{-1}+\ldots)$. En faisant la division suivant les puissances d\'ecroissantes, on a 
$$
\frac{d\widetilde{F}}{\widetilde{F}} \wedge \omega = \frac{da_k}{a_k}\wedge \omega +k\frac{dG\wedge \omega}{G}+\frac{d(a_1 G^{-1}+\ldots)}{(1+a_1 G^{-1}\ldots)}\wedge \omega.
$$
Le dernier terme de la somme contient des puissances de $G$ inf\'erieures ou
\'egales \`a $-1$. Le deuxi\`eme terme est de degr\'e $-1$ en $G$. Par
transcendance de $G$ on en d\'eduit que $\frac{da_k}{a_k}\wedge \omega = p d
\omega$. A partir du facteur int\'egrant $\sqrt[p]{a_k}$ on construit une
suite d'extensions de longueur $n-1$ contenant une int\'egrale
premi\`ere. Ceci contredit \`a nouveau la minimalit\'e de la suite.\\
\indent On a obtenu une suite de longueur deux $K_{0} \subset K_1 \subset
K_2$ avec $K_1=K_0(\sqrt[p]{\widetilde{F}})$ et $K_2=K_1(H)$ o\`u
$dH=\sqrt[p]{\widetilde{F}}\omega$ : il existe un élément de $K$ dont
une racine est un facteur int\'egrant pour $\omega$.

R\'eciproquement, une suite de Godbillon-Vey de longueur un permet par ces
formules de construire une int\'egrale premi\`ere de type Darboux pour $\omega$.

 Ceci achève la peuve du théorème.
\end{preuve}

%%%%%%%%%%%%%%%%%%%%%%%%%%%%%%%%%%%%%%%%%%%%%%%%%%%%%%%%%%%%%%%%%%%%%%%%%%%%%%%%%%%%%%%%%%%%%%%%%%%%%%
%%%
\section{Le groupo\"ide de Galois d'un germe de feuilletage de $(\mathbb{C}^{2},0)$ \`{a} singularit\'{e} r\'{e}duite}

Nous allons maintenant d\'eterminer les feuilletages sur un bidisque \`{a}
singularit\'{e} r\'{e}duite dont le groupo\"ide de Galois est de rang
transverse fini en fonction des invariants de leurs classes analytiques. 
\begin{definition}
Un feuilletage $\mathcal{F}$ de $(\mathbb{C}^{2},0)$ sera dit \`{a}
singularit\'{e} r\'{e}duite si il existe une forme $\omega$ d\'efinissant $\mathcal{F}$ dont la partie lin\'{e}aire s'\'{e}crit dans de bonnes coordonn\'{e}es:
  \begin{enumerate}
   \item $\lambda_{1}xdy+\lambda_{2}ydx ,$ $(\lambda_{1},\lambda_{2}) \in \mathbb{C}^{\ast} \times \mathbb{C}^{\ast}$, $\frac{\lambda_{1}}{\lambda_{2}} \notin \mathbb{Q}_{<0}$,
   \item $ydx$. 
  \end{enumerate}
\end{definition}
La terminologie employ\'ee renvoie au th\'eor\`eme de r\'eduction de
Seidenberg \cite{seidenb} : ces singularit\'es sont les plus simples que
l'on obtient apr\`es \'eclatements.\\
%Les r\'ef\'erences que nous utilisons sont \cite{Du}, \cite{MM}, \cite{MR1}, \cite{MR2}.\\
\indent Les feuilletages de type (1) sont appel\'es des selles. Ils admettent deux
courbes analytiques invariantes lisses et transverses dont les holonomies ont
pour parties lin\'eaires $e^{-2i\pi \lambda_{1}/ \lambda_{2}}$ et $e^{-2i\pi
  \lambda_{2}/ \lambda_{1}}$. Ces feuilletages ont des comportements
diff\'erents suivant les valeurs de $\frac{\lambda_1}{\lambda_2}$. Lorsque
$\frac{\lambda_1}{\lambda_2}$ n'est pas r\'eel ou r\'eel n\'egatif non
rationnel, on sait, d'apr\`es Poincar\'e, que le feuilletage est analytiquement lin\'earisable. Lorsque 
$\frac{\lambda_1}{\lambda_2}$ est r\'eel irrationnel, le feuilletage est
formellement lin\'earisable \cite{Ilyashenko-dansMR1}. Lorsque
$\frac{\lambda_1}{\lambda_2}$ est rationnel ces feuilletages s'appellent selles
r\'esonnantes. Ils ne sont plus lin\'earisables, mais admettent les formes normales formelles suivantes :
$$ 
p(1+(\lambda -1)(x^{p}y^{q})^{k})ydx + q(1+\lambda(x^{p}y^{q})^{k})xdy.
$$
Les axes de coordonn\'ees sont des courbes invariantes pour ces feuilletages. 
L'holonomie de la feuille $x=0$ calcul\'ee sur la transverse $(y=1,x=t)$ est le diff\'eomorphisme 
$$
f=e^{-2i\pi p/q}\exp\left(\frac{t^{qk+1}}{1-\lambda\frac{q}{p} t^{qk}}\frac{d}{dt}\right).
$$
% et celle de la feuille $x=0$ est $e^{-2i\pi
% q/p}\exp(\frac{x^{pk+1}}{1-(\lambda-1)\frac{p}{q} x^{pk}}\frac{d}{dx})$.
La classe de conjugaison analytique de cette holonomie est un invariant
complet de la classe analytique du feuilletage.
D'apr\`es les r\'esultats de \cite{MR2}, apr\`es pr\'eparation du feuilletage,
il existe des normalisations analytiques $h_i$ sur des secteurs $U_i$ de la forme
$-\frac{\pi}{2}- \varepsilon < \arg((x^py^q)^k) < \frac{\pi}{2}+ \varepsilon$
qui sont asymptotes \`a la normalisante formelle tangente \`a l'identit\'e. Le
cocycle $(U_i \cap U_{i+1}, h_{i+1} \circ h_i^{-1})_{i}$ induit un invariant
complet de la classe analytique du feuilletage.\\
\indent Les feuilletages de type (2) sont appel\'es des n{\oe}ud-cols. Ils ont pour formes normales formelles 
$$ 
x^{k+1}dy-y(1-\lambda x^{k})dx 
$$
d'holonomie $f=\exp(\frac{t^{k+1}}{1-\lambda t^{k}}\frac{d}{dt})$ calcul\'ee sur la transverse $(y=1,x=t)$ \`a la s\'eparatrice forte $x=0$.
D'apr\`es les r\'esultats de \cite{MR1}, apr\`es pr\'eparation du feuilletage,
il existe des normalisations analytiques $h_i$ sur des secteurs $U_i$ de la forme $-\frac{\pi}{2}- \varepsilon < \arg(x^k) < \frac{\pi}{2}+
\varepsilon$ qui sont asymptotes \`a la normalisante formelle tangente \`a
l'identit\'e. Elles d\'efinissent un cocycle qui est un invariant complet de
la classe analytique du feuilletage. \\
Dans le cas des selles r\'esonnantes et des n{\oe}ud-cols, les invariants de la classe de conjugaison analytique de l'holonomie de la s\'eparatrice forte sont donn\'es par les composantes transverses des invariants analytiques du feuilletage.

\begin{proposition}
\label{extension}
Soit $f$ l'holonomie d'un feuilletage \`a singularit\'e r\'eduite
$\mathcal{F}$. Si $f$ est solution d'un \dgl\ $\G$ sur un disque
transverse $T$ alors il existe un \dgl\ admissible pour $\mathcal{F}$ dont
l'\'equation transverse au voisinage de $T$ co\"incide avec celle de $\G$.  
\end{proposition}

\begin{preuve}
Nous commen\c{c}ons par d\'eterminer la liste des \dgls\ de rang transverse sup\'erieur ou
\'egal \`a un admissible pour le n{\oe}ud-col mod\`ele  $\omega = x^{2}dy-ydx$.  
Le facteur int\'egrant $F = \frac{1}{x^2y}$ d\'etermine un \dgl\ de rang transverse un
obtenu en \'ecrivant l'invariance de la forme ferm\'ee $F \omega$ :
$$
\frac{dY}{Y}-\frac{dX}{X^2}=\frac{dy}{y}-\frac{dx}{x^2}.
$$
Consid\'erons les cartes $y \not \in \RR_{\leq 0}$ et $y \not \in \RR_{\geq
  0}$. En choisissant deux d\'eterminations de $\log y$, on les munit des
coordonn\'{e}es redressantes $t=\frac{x}{x \log y+1}$ et $z=y$. Dans ces
coordonn\'ees, on v\'erifie que les \'equations de ce \dgl\ s'\'ecrivent :  
$$
\frac{\partial T}{\partial z}=0 \text{ et } \frac{1}{T^{2}} \frac{\partial T}{\partial t} = \frac{1}{t^{2}}.
$$
% These detail
\noindent On retrouve la formule (1) de la remarque \ref{formules} qui donne sur la partie transverse le \dgl\ $\G_1(-\frac{2}{t})$.

Les \dgls\ de rang transverse deux admissibles pour $\omega$ sont obtenus d'apr\`es le th\'eor\`eme
\ref{suitegv} \`a partir de toutes les formes ferm\'ees $\alpha$ v\'erifiant $d\omega =
\omega \wedge \alpha$. Dans le cas du n{\oe}ud-col, ces formes s'\'ecrivent  
$$
\alpha = \frac{dF}{F} + cF\omega
$$
o\`u $c$ est
un nombre complexe quelconque.   
Consid\'erons les cartes $y \not \in \RR_{\leq 0}$ et $y \not \in \RR_{\geq
  0}$.  Dans les
coordonn\'ees pr\'ec\'edentes, $\omega=w(z,t)dt$.
D'apr\`es la formule (2) de la remarque \ref{formules}, l'\'equation transverse du \dgl\ associ\'e \`a la suite de Godbillon-Vey $(\omega,\alpha)$ est $\G_2(\mu)$ avec $\mu(t)dt= \alpha + \frac{dw}{w}$.
On v\'erifie que les \'equations de ce \dgl\ s'\'ecrivent :  
$$
\frac{\partial T}{\partial z}=0 \text{ et } \G_2(-\frac{c}{t^2}-\frac{2}{t}).
$$
% Le \dgl\ $\G_2(\frac{c}{t^2}-\frac{2}{t})$ contient $\G_1(-\frac{2}{t})$.

Les \dgls\ de rang transverse trois sont obtenus en prenant toutes les formes
$\alpha$ v\'erifiant la premi\`ere \'equation de Godbillon-Vey et en
compl\'etant la suite par l'unique forme $\beta$ satisfaisant les deux
derni\`eres \'equations de Godbillon-Vey. On obtient les suites $(\omega =
x^{2}dy-ydx, \alpha =\frac{dF}{F} + cF\omega, \beta=F^2\omega)$. Dans les
cartes de coordonn\'ees $(t,z)$ pr\'ec\'edentes, d'apr\`es la formule (3) de la remarque \ref{formules}, l'\'equation transverse du \dgl\ d\'efinie par la suite $(\omega,\alpha,\beta)$ est $\G_3(\nu)$ avec 
$$
\nu(t)dt = w \beta + d\left(\frac{\partial_t w}{w}+a_t\right)-\frac{1}{2}\left(\frac{\partial_t w}{w}+a_t\right)^2dt 
$$   
o\`u $\omega = w(z,t)dt$ et $\alpha = a_tdt+a_zdz$.
On en d\'eduit que les \'equations de ce \dgl\ s'\'ecrivent :  
$$
\frac{\partial T}{\partial z}=0 \text{ et } \G_3(-\frac{c^2}{t^4}).
$$

D'autre part, d'apr\`es le d\'ebut de la section \ref{sectiondiffeo}, les \dgls\ sur le disque transverse contenant l'holonomie $\exp(t^2\frac{d}{dt})$
du n{\oe}ud-col sont  $\G_1(-\frac{2}{t})$, $\G_2(-\frac{c}{t^2}-\frac{2}{t})$
et $\G_3(-\frac{c^2}{t^4})$ et seulement ceux-ci. Ces deux listes \'etant
identiques, nous avons montr\'e la proposition pour le n{\oe}ud-col $x^2dy-ydx$.
Les autres formes normales de n{\oe}ud-cols, $x^{k+1}dy-y(1-\lambda x^{k})dx$,
se ram\`enent au cas pr\'{e}c\'{e}dent par $(x,y) \mapsto (\frac{x}{1-\lambda
  x \log x},y)$ et la ramification $(x,y) \mapsto (x^{k},y)$. Les formes
normales de selles r\'{e}sonnantes se ram\`{e}nent aux n{\oe}ud-cols par
l'\'{e}clatement $(x,y) \mapsto (xy,y)$ et la ramification $(x,y) \mapsto
(x^{p},y^{q})$. Tous ces feuilletages admettant des facteurs int\'egrants, on
peut aussi d\'eriver directement les \'equations des \dgls\ admissibles pour un
de ces feuilletages en consid\'erant toutes les suites de Godbillon-Vey que
l'on peut associer \`a ce feuilletage, \`a \'equivalence pr\`es.\\

Consid\'erons un n{\oe}ud-col dans la classe formelle de $x^{k+1}dy-y(1-\lambda x^{k})dx$ dont
l'holonomie n'est plus analytiquement normalisable, et supposons que celle-ci
soit solution d'un \dgl\ de rang trois $\G_3(\nu)$. Nous allons construire un \dgl\ admissible pour le feuilletage de rang transverse trois. 
%La fin de la preuve utilisera les constructions faites dans \cite{MR1} et
%\cite{MR2} ainsi que dans \cite{C2}.
On sait d'apr\`es le th\'eor\`eme \ref{5}, que $f$ est
solution d'un \dgl\ de rang trois si et seulement si sa forme normale formelle et
ses invariants analytiques sont eux-m\^emes solutions d'un \dgl\ de rang trois $\G_3(\overline{\nu})$. Soit $\widehat{h}$ la conjugante formelle entre $f$ et sa forme
normale. On a
$$\nu=\overline{\nu}\circ\widehat{h}(\widehat{h}')^2+S\widehat{h}.$$
Nous avons prouv\'e ci-dessus que ce \dgl\ se prolonge en un \dgl\ $\G$
admissible pour le feuilletage mod\`ele de rang transverse trois.

Maintenant nous allons construire un \dgl\ admissible pour le feuilletage
initial \`a partir de $\G$. Quitte \`a faire une conjugaison analytique, on peut
supposer que la s\'eparatrice forte du feuilletage a pour \'equation $x=0$. Le
n{\oe}ud-col est alors conjugu\'{e} au-dessus des secteurs
$(-\frac{\pi}{2}-\epsilon \le \arg(x^{2k}) \le \frac{\pi}{2}+\epsilon)$ \`{a} sa
forme normale formelle par des normalisantes sectorielles $h_i$, avec $i$ variant dans $\ZZ / 2k\ZZ$,
asymptotes \`a la normalisante formelle. On consid\`ere alors les \dgls\
$h_i^*\G$ au-dessus de chaque secteur. Une fois que l'on s'est fix\'e les deux premières formes $(\omega,\alpha)$ d'une suite de Godbillon-Vey de logueur trois pour le n{\oe}ud-col, ce \dgl\ est la donn\'ee d'une troisième forme $\beta_i$ satisfaisant les \'equations de Godbillon-Vey.
Montrons que sur les intersections de deux de ces secteurs les deux formes $\beta_i$ et $\beta_{i+1}$ co\"incident. Puisque, d'apr\`es \cite{MR1}, les composantes transverses de $h_i \circ
h_{i+1}^{-1}$ sont les composantes du cocycle des invariants de l'holonomie $f$, celui-ci \'etant solution de $\G_3(\overline{\nu})$, l'automorphisme du
feuilletage mod\`ele $h_i \circ h_{i+1}^{-1}$ est solution de $\G$.
Soit $(t,z)$ des coordonn\'ees redressantes au voisinage d'un disque $(z=0)$ transverse \`a la s\'eparatrice forte $(t=0)$. En \'ecrivant les \'equations des \dgls\ correspondants aux triplets $(\omega,\alpha, \beta_i)$ et $(\omega,\alpha, \beta_{i+1})$ (formules \ref{formules}) on obtient pour chacune des \'equations sur les secteurs transverses correspondant,
$
\G_3(\nu_i) \text{ et }\G_3(\nu_{i+1})
$
avec 
$$\nu_i = \overline{\nu}\circ h_i ( h_i')^2+S(h_i)$$
où on d\'esigne par $h_i$ la composante transverse de la normalisante sectorielle $h_i$. 
Comme la composante transverse de $h_i \circ h_{i+1}^{-1}$ est solution de $\G_3(\overline{\nu})$ $\nu_i=\nu_{i+1}$, les fonctions $\nu_i$ et $\nu_{i+1}$ \'etant asymptotes \`a $\nu$, elles sont \'egales \`a cette derni\`ere. La forme $\beta$ est ainsi bien définie et m\'eromorphe.

%% METTRE LES INDICES AU-DESSUS

%% Les autres n{\oe}ud-col se traitent de la m\^eme mani\`ere. Apr\`es pr\'eparation, un n{\oe}ud-col est conjugu\'{e} au-dessus de secteur $(-\frac{\pi}{2}-\epsilon \le \arg(x^{k}) \le \frac{\pi}{2}+\epsilon)$ \`{a} sa forme normale formelle. Nous noterons $h_i$ la normalisante sur le $i$-i\`eme secteur et $\G$ le \dgl\ admissible pour la forme normale contenant les invariants analytiques du n{\oe}ud-col.
%% On ram\`ene les \'equations de ce $\mathcal{D}$-groupo\"ide de Lie sur
%% chaque secteur de conjuguaison ce qui nous donne $2k$ \dgls\, $h_{i}^* \G$
%% au-dessus des $2k$ secteurs de conjugaisons tous admissibles pour le
%% n{\oe}ud-col restreint au secteur. Comme les invariants analytiques du
%% feuilletage sont solutions du $\mathcal{D}$-groupo\"ide de Lie mod\`ele, on
%% a $(h_i^{-1})^*\ h_{i-1}^* \G=\G$ donc $h_i^* \G = \ h_{i-1}^* \G$. Les
%% \'equations sectorielles obtenues se recollent en des \'equations au-dessus
%% du bidisque priv\'e de l'axe $\{x=0\}$. Les conjugantes admettant un
%% d\'eveloppement asymptotique le long de l'axe des $y$, nous obtenons des
%% \'equations  m\'eromorphes par le th\'eor\`eme des singularit\'es
%% inexistantes.\\

 Les selles r\'esonnantes se traitent exactement de la m\^{e}me mani\`{e}re
: seuls les secteurs changent de formes et sont donn\'es par $(-\frac{\pi}{2}-\epsilon \le \arg((x^{p}y^{q})^k) \le \frac{\pi}{2}+\epsilon)$. 

La preuve dans le cas d'une holonomie unitaire est analogue.     
\end{preuve}

\begin{proposition}
Soit $\mathcal{F}$ un germe de feuilletage \`a singularit\'e r\'eduite.  
\begin{enumerate}
\item $\mathcal{F}$ admet un $\mathcal{D}$-groupo\"ide de Lie admissible de
  rang transverse un si et seulement son holonomie est analytiquement normalisable.
\item $\mathcal{F}$ admet un $\mathcal{D}$-groupo\"ide de Lie admissible de rang transverse deux si et seulement si son holonomie est unitaire.
\item $\mathcal{F}$ admet un $\mathcal{D}$-groupo\"ide de Lie admissible de rang transverse trois si et seulement si son holonomie est binaire.
\end{enumerate}
\end{proposition}
\begin{preuve}
Soit $\mathcal{F}$ un feuilletage r\'eduit admettant un \dgl\ de rang transverse
fini. Celui-ci d\'efinit un \dgl\ de m\^eme rang contenant l'holonomie par le
lemme \ref{local}. Lorque l'holonomie est un
diff\'eomorphisme r\'esonnant, le
th\'eor\`eme \ref{theoimp} nous assure qu'elle est normalisable,
unitaire ou binaire suivant la valeur du rang. Dans le cas des holonomies formellement
lin\'earisables, le
th\'eor\`eme \ref{pasdepetitdiv} nous assure qu'elle est analytiquement
lin\'earisable. %Un th\'eor\`eme de Mattei-Moussu \cite{MM} assure que dans ce
                %cas le feuilletage est analytiquement lin\'earisable. 
%% Pour la r\'eciproque, consid\'erons d'abord les feuilletages formellement
%% lin\'earisables qui admettent un \dgl\ de rang transverse fini. D'apr\`es le
%% th\'eor\`eme \ref{pasdepetitdiv} son holonomie est analytiquement
%% lin\'earisable. Un th\'eor\`eme de Mattei-Moussu \cite{MM} assure que dans ce
%% cas le feuilletage est analytiquement lin\'earisable. 
Pour la r\'eciproque, consid\'erons d'abord les feuilletages analytiquement lin\'earisables
ou normalisables. Il admettent toujours un facteur int\'egrant et donc un
\dgl\ admissible de rang transverse un.

Consid\'erons ensuite les selles r\'esonnantes et les n{\oe}ud-cols
d'holonomie unitaire ou binaire. Le \dgl\ donn\'e par le th\'eor\`eme
\ref{theoimp} s'\'etend gr\^ace \`a la proposition \ref{extension}
pr\'ec\'edente en un \dgl\ admissible pour le feuilletage de rang transverse deux ou trois. 
\end{preuve}
Nous obtenons ainsi une nouvelle preuve de la caract\'erisation sur les invariants analytiques des feuilletages \`a singularit\'e r\'eduite admettant une structure transverse affine m\'eromorphe ou une structure transverse projective m\'eromorphe.    

Ces r\'esultats ont d\'ej\`a \'et\'e obtenus en utilisant d'autres techniques dans \cite{BT} pour le cas transversalement affine et \cite{T} dans le cas transversalement projectif.
% maintenir ??

%%%%%%%%%%%%%%%%%%%%%%%%%%%%%%%%%%%%%%%%%%%%%%%%%%%%%%%%%%%%%%%%%%%%%%%%%%%%%%%%%%%%%%%%%%%%%%%%%%%%
%%      Groupo\"ides de Galois et extensions fortement normales}
\section{Groupo\"ides de Galois et extensions fortement normales}
\label{extensionfn}

Dans cette section $\mathcal{M}_\Delta$ d\'esigne le corps des fonctions m\'eromorphes sur le polydisque $\Delta$ de $\CC^n$.
Consid\'erons l'espace $J_k^*(\Delta \to \mathbb{C})$ des jets d'ordre $k$
d'applications submersives de $\Delta$ dans $\mathbb{C}$ (la notation $*$ d\'esigne ici la propri\'et\'e de submersivit\'e). Le choix de
coordonn\'ees $x$ sur $\Delta$ et $H$ sur $\mathbb{C}$ nous permet
d'identifier cet espace \`a un ouvert de $\Delta \times \mathbb{C}
\times_{|\alpha|\leq k}\mathbb{C}^\alpha$ avec les coordonn\'ees $H^\alpha$
naturellement associ\'ees au choix de $x$ et de $H$. Ces espaces sont munis de l'anneau des \'equations aux d\'eriv\'ees partielles d'ordre inf\'erieur \`a $k$, $\mathcal{O}_{J_k^*(\Delta \to \mathbb{C})}=\mathcal{O}_{\Delta}[H,\ldots,H^\alpha\ldots]$ et pour chaque d\'erivation partielle $\frac{\partial}{\partial x_i}$ d'une d\'erivation $D_i : \mathcal{O}_{J_k^*(\Delta \to \mathbb{C})} \to \mathcal{O}_{J_{k+1}^*(\Delta \to \mathbb{C})}$.

\begin{definition}
Une $\mathcal{D}$-vari\'et\'e dans $J^*(\Delta \to \mathbb{C})=\underset{\longleftarrow}{lim}\ J_k^*(\Delta \to \mathbb{C}) $ est donn\'ee par un id\'eal $\mathcal{J} \subset \mathcal{O}_{J^*(\Delta \to \mathbb{C})} = \underset{\longrightarrow}{lim}\ \mathcal{O}_{J_k^*(\Delta \to \mathbb{C})}$ diff\'erentiel et r\'eduit tel que $\mathcal{J} \cap \mathcal{O}_{\Delta} =\emptyset$.
\end{definition}

\begin{definition}
Soit $\mathcal{A}$ un anneau diff\'erentiel sur $\mathcal{M}_\Delta$. Le
spectre diff\'erentiel est l'ensemble $Spec^{diff} (\mathcal{A})$ des id\'eaux premiers et diff\'erentiels de $\mathcal{A}$.
\end{definition}
Cet ensemble peut \^etre muni d'une topologie appel\'ee topologie de
Zariski-Kolchin (\cite{B}, \cite{K}).
Nous allons \'etudier dans cette partie la $\mathcal{D}$-vari\'et\'e des
int\'egrales premi\`eres d'un germe de feuilletage $\mathcal{F}$ d\'efini par
des formes $\omega_1, \ldots \omega_q$. Elle est
donn\'ee par le syst\`eme d'\'equations aux d\'eriv\'ees partielles $dH_i \wedge
\omega_1\wedge\ldots\wedge\omega_q =0 $. Consid\'erons l'anneau diff\'erentiel 
$$
\M_{\Delta}\otimes_{\O_{\Delta}}\O_{J^*(\Delta \to \mathbb{C}^q)}/(dH_i \wedge
\omega_1\wedge\ldots\wedge\omega_q)
$$
o\`u $(dH_i \wedge \omega_1\wedge\ldots\wedge\omega_q)$ est l'id\'eal diff\'erentiel r\'eduit engendr\'e par les
composantes de ces $(q+1)$-formes pour $i$ entre $1$ et $q$, 
et notons $\mathcal{O}_{\F}$ son localis\'e sur $dH_1\wedge\ldots
dH_q \not = 0$.

L'ensemble $Spec^{diff}\mathcal{O}_{\F}$ repr\'esente l'ensemble des
syst\`emes d'\'equations aux d\'eriv\'ees partielles compatibles avec le fait
d'\^etre un syst\`eme complet d'int\'egrales premi\`eres. La notion de r\'eductibilit\'e d'un syst\`eme
d'\'equations aux d\'eriv\'ees partielles de Jules Drach (\cite{Dr}) correspond
\`a la non trivialit\'e du spectre diff\'erentiel. 
\begin{definition}[\cite{buium},\cite{Ko}, \cite{K}]
Soit $\mathcal{M}_\Delta \subset \mathcal{K}$ une extension de degr\'e de
transcendance fini de corps
diff\'erentiels. Cette extension sera dite fortement normale si pour toute
extension diff\'erentielle $\mathcal{E}$ de $\mathcal{K}$ et tout morphisme $\sigma : \mathcal{K} \to \mathcal{E}$ au-dessus de $\mathcal{M}_\Delta$ :
\begin{enumerate}
\item $\sigma$ laisse les constantes de $\mathcal{K}$ invariantes,
\item $\sigma (\mathcal{K}) \cdot \mathcal{E}^c = \mathcal{K} \cdot \mathcal{E}^c$ 
\end{enumerate}
o\`u $\mathcal{E}^c$ d\'esigne le corps des constantes de $\mathcal{E}$ et le point d\'esigne le compositum des corps dans $\mathcal{E}$.
\end{definition}
D'apr\`es la th\'eorie de Kolchin (voir \cite{Ko} et \cite{K}), le groupe de Galois de ces extensions est un groupe alg\'ebrique. Ses
sous-groupes alg\'ebriques sont en correspondance avec les extensions
diff\'erentielles interm\'ediaires. 

Le th\'eor\`eme suivant confirme les r\'esutats incomplets de J. Drach \cite{Dr2} et s'inscrit dans ``une th\'eorie g\'en\'erale de la r\'eductibilit\'e des \'equations'' esquiss\'ee par E. Vessiot \cite{Ve1}, \cite{Ve2}.

\begin{theoreme}
\label{extfortnorm}
Soit $\mathcal{F}_\omega$ un feuilletage de codimension un de $(\mathbb{C}^n,0)$. Les assertions suivantes sont \'equivalentes :
\begin{enumerate}
\item[\emph{(1)}] le groupo\"ide de Galois de $\mathcal{F}_\omega$ est propre;
\item[\emph{(2)}] le spectre diff\'erentiel de $\mathcal{O}_{\F}$ est non trivial : $Spec^{diff} (\mathcal{O}_{\F}) \not = \{0\}$;
\item[\emph{(3)}] il existe une int\'egrale premi\`ere de $\F$ dans une extension fortement normale $\mathcal{K}$ de $\mathcal{M}_\Delta$. 
\end{enumerate}
\end{theoreme}

Nous d\'emontrerons successivement les implications (3) $\Rightarrow$ (2), (2)
$\Rightarrow$ (1) et (1) $\Rightarrow$ (3) dans les lemmes suivants. Nous
montrerons les deux premi\`eres pour un feuilletage de codimension quelconque.

\begin{lemme}
Soit $\F$ un feuilletage donn\'e par $q$ $1$-formes. Si il existe $q$ int\'egrales
premi\`eres fonctionnellement ind\'ependantes dans une extension
diff\'erentielle de $\M_{\Delta}$ de degr\'e de transcendance fini alors le spectre diff\'erentiel de $\mathcal{O}_{\F}$ est non trivial.  
\end{lemme}
\begin{preuve}
L'existence d'un syst\`eme d'int\'egrales premi\`eres dans $\mathcal{K}$ donne un morphisme
diff\'erentiel au-dessus de $\mathcal{M}_\Delta$ :
$$
\mathcal{O}_{\F} \longrightarrow \mathcal{K}
$$
induit par l'identification des coordonn\'ees $H_i$ avec
les int\'egrales premi\`eres.
Le noyau de ce morphisme est un id\'eal diff\'erentiel premier de
$\mathcal{O}_{\F}$ et donne donc un \'el\'ement de
$Spec^{diff} (\mathcal{O}_{\F})$. L'extension $\mathcal{K}$ \'etant de degr\'e de transcendance fini ce qui n'est pas le
cas de $\mathcal{O}_{\F}$, le
morphisme ne peut pas \^etre injectif. L'\'el\'ement obtenu dans
$Spec^{diff}(\mathcal{O}_{\F})$ est non trivial. 
\end{preuve}

\begin{lemme}
\label{reduction}
Soit $\F$ un feuilletage donn\'e par $q$ $1$-formes. Si le spectre
diff\'erentiel de $\mathcal{O}_{\F}$ est non trivial alors le groupo\"ide de Galois de $\mathcal{F}_\omega$ est propre.
\end{lemme}
\begin{preuve}
Soit $\J$ un id\'eal diff\'erentiel premier de $\O_{J^*(\Delta \to \CC^q)}$
contenant l'id\'eal diff\'erentiel donn\'e par $dH_i \wedge
\omega_1 \wedge\ldots \omega_q$. Nous
allons construire un \dgl\ dont les solutions sont les germes $\varphi$ tels
que pour tout $H=(H_1,\ldots, H_q)$, $H$ est solution de $\J$ si et seulement
  si $H \circ \varphi$  est solution de $\J$.
Pour cela on consid\`ere l'action de $J_p^*(\Delta)$ sur $J_p^*(\Delta \to
\mathbb{C}^q)$ par composition \`a la source :
$$
comp : J^*(\Delta \to \mathbb{C}^q)\times_{\Delta}J^*(\Delta) \longrightarrow J^*(\Delta \to \mathbb{C}^q).
$$ 
Cette action se traduit sur les anneaux par l'existence d'une fl\`eche $comp^*$
 satisfaisant les diagrammes commutatifs suivant :
\begin{itemize}
  \item l'identit\'e
\vspace*{-1.3cm}
 $$ 
  \shorthandoff{;:!?}
  \xymatrix{
\relax  \\
\mathcal{O}_{J^*(\Delta \to \mathbb{C}^q)} \ar[rr]^-{comp^*} \ar[rrd]^{=}  &&   \mathcal{O}_{J^*(\Delta \to \mathbb{C}^q)} \otimes_{\mathcal{O}_{\Delta}} \mathcal{O}_{J^*(\Delta)} \ar[d]^{1 \otimes e^*}  \\
    &&  \mathcal{O}_{J^*(\Delta \to \mathbb{C}^q)}
}
  $$
  \item la composition
\vspace*{-1.3cm}
$$ 
  \shorthandoff{;:!?}
  \xymatrix{
\relax  \\
\mathcal{O}_{J^*(\Delta \to \mathbb{C}^q)} \ar[rr]^-{comp^*} \ar[d]_{comp^*}  &&  \mathcal{O}_{J^*(\Delta \to \mathbb{C}^q)} \otimes_{\mathcal{O}_{\Delta}} \mathcal{O}_{J^*(\Delta)} \ar[d]^{1 \otimes c^*}  \\
\mathcal{O}_{J^*(\Delta \to \mathbb{C}^q)} \otimes_{\mathcal{O}_{\Delta}} \mathcal{O}_{J^*(\Delta)} \ar[rr]^-{comp^* \otimes 1}  &&  \mathcal{O}_{J^*(\Delta \to \mathbb{C}^q)} \otimes_{\mathcal{O}_{\Delta}} \mathcal{O}_{J^*(\Delta)} \otimes_{\mathcal{O}_{\Delta}} \mathcal{O}_{J^*(\Delta)}
}
  $$

\end{itemize}
Un germe de diff\'eomorphisme $\varphi : (\Delta,a) \to (\Delta,b)$ induit par prolongement un morphisme $\varphi^*$ des jets de $J^*(\Delta \to \mathbb{C}^q)$ de source $b$
sur ceux de source $a$ par composition. \'Etant donn\'e un id\'eal
$\mathcal{J}$ de $\mathcal{O}_{J_p^*(\Delta \to \mathbb{C}^q)}$, nous allons
chercher \`a d\'eterminer les \'equations diff\'erentielles satisfaites par
les germes $\varphi$ tels que $\varphi^{**} (\mathcal{J}\otimes \mathbb{C}_a) =
(\mathcal{J}\otimes \mathbb{C}_b)$ (on note $\varphi^{**}$ le morphisme d'anneau induit par la transformation $\varphi^*$). Consid\'erons l'id\'eal $comp^* \J$ dont
les solutions sont l'ensemble des couples $(H,\varphi)$ tels que $H \circ
\varphi$ est solution de $\J$. Les solutions de $comp^*\mathcal{J} +
\mathcal{J}\otimes 1$ sont les couples $(H,\varphi)$ tels que $H$ et $H \circ
\varphi$ soient solutions de $\J$. Il faut d\'eterminer le plus petit id\'eal
$\mathcal{I}$ de $\mathcal{O}_{J^*(\Delta)}$ v\'erifiant $comp^*\mathcal{J} +
\mathcal{J}\otimes 1 \subset \mathcal{J} \otimes 1 + 1 \otimes \mathcal{I}$ en
dehors d'une hypersurface. Les solutions de $\I$ sont les $\varphi$ tels que
si $H$ est solution de $\J$ alors $H \circ \varphi$  est solution de $\J$. \\
Nous allons d\'eterminer un syst\`eme de g\'en\'erateurs de l'id\'eal $\I$.
En l'absence de torsion, nous noterons encore $\mathcal{J}$ l'id\'eal engendr\'e par $\mathcal{J}$ dans
$\M_{\Delta} \otimes \O_{J^*(\Delta \to \CC^q)}$. D'apr\`es le th\'eor\`eme de ``noetherianité'' \ref{malg-ritt-rad}, il existe un entier $k$ tel que
l'id\'eal $\mathcal{J}$ soit
diff\'erentiablement engendr\'e, en tant qu'id\'eal r\'eduit, par ses \'el\'ements d'ordre inf\'erieur \`a $p$. On note
$\mathcal{J}_p$ la trace de l'id\'eal $\J$ dans l'anneau des \'equations
diff\'erentielles d'ordre inf\'erieur ou \'egale \`a $p$ : $\O_{J_p^*(\Delta \to \CC^q)}$.\\
Soit $f_1, \ldots ,f_n$ un syst\`eme g\'en\'erateur de l'id\'eal
$\mathcal{J}_p$, qu'on supposera $\mathcal{M}_\Delta$-libre. On le compl\`ete
en une base du $\mathcal{M}_\Delta$-espace vectoriel $\mathcal{J}_p$ : $f_1,
\ldots ,f_n, $ $\ldots, f_m, \ldots$ Puis on compl\`ete cette famille par des $e_1,\ldots, e_k,
\ldots$ en une $\mathcal{M}_\Delta$-base de $\M_{\Delta} \otimes \O_{J_p^*(\Delta
  \to \CC^q)}$. Soit $f \in \M_{\Delta} \otimes \O_{J_p^*(\Delta \to \CC^q)} $, nous noterons :
$$
comp^*f=\sum f_j\alpha^j(f) + \sum e_k\beta^k(f),
$$
avec les $\alpha^j(f)$ et $\beta^k(f)$ dans $\M_\Delta \otimes _{\O_{\Delta}}\mathcal{O}_{J_p^*(\Delta)}$.
On consid\`ere alors l'id\'eal $ \mathcal{I}_p$ de $\M_\Delta \otimes _{\O_{\Delta}}\mathcal{O}_{J_k^*(\Delta)}$ d\'efini par les $\beta^k(f)$ pour $f\in \mathcal{J}_p$. Par construction cet id\'eal v\'erifie deux propri\'et\'es importantes.

\begin{enumerate}

\item[(a)] Il est engendr\'e par les $\beta^k(f_i)$ pour $i=1, \ldots ,n$. En effet,
    on a clairement $\beta^k(f+g) = \beta^k(f)+ \beta^k(g)$. D'autre part
    comme $$ comp^*(fg) \equiv \sum e_\ell e_k \beta^\ell(f)\beta^k(g) \mod
    \J_p\otimes 1$$ en \'ecrivant $e_\ell e_k$ dans la base d\'ecrite au-dessus,
    on obtient que $\beta^i(fg)$ est une combinaison des $\beta^\ell(f)\beta^k(g)$
    \`a coefficients dans $\M_\Delta$. En particulier, quelque soit $f$ dans
    $\J_p$, $\beta^k(f)$ est dans l'id\'eal engendr\'e par les $\beta^k(f_i)$
    pour $i=1, \ldots ,n$. 

\item[(b)] En consid\'erant la d\'ecomposition de $comp^* f_i$ on remarque que, pour
tout $\varphi$ de $J_p^*(\Delta)$ solution de $\I_p$, si $H$ est solution de
$\J_p$ alors $H\circ \varphi$ est aussi solution de $\J_p$. R\'eciproquement
si pour toute solution $H$ de $\J_p$,  $H\circ \varphi$ est encore solution de
$\J_p$ alors, en utilisant l'ind\'ependance des $e_k$, on a $\varphi$ solution
de $\I_p$. Un jet $\varphi$ est solution de $\I_p$ si et seulement si
$\varphi^{**} \mathcal{J}_p\otimes \mathbb{C}_a \subset \mathcal{J}_p\otimes
\mathbb{C}_b$.
\end{enumerate} 

V\'erifions maintenant que cet id\'eal v\'erifie les propri\'et\'es (1) et (3)
d'un \dgl. L'inclusion de cet id\'eal dans l'id\'eal d\'efinissant
l'identit\'e se d\'eduit du diagramme commutatif de l'identit\'e
ci-dessus. D'autre part on a :
\begin{multline}
(comp^* \otimes 1)(comp^*)(f)= 
  \sum _{i,j}f_i\alpha^i(f_j)\alpha^j(f)  \sum _{l,j}
   e_l\beta^l(f_j)\alpha^j(f) \nonumber \\  + \sum _{i,k}f_i\alpha^i(e_k)\beta^k(f) +
   \sum _{k,l}e_l\beta^l(e_k)\beta^k(f)\nonumber
\end{multline} 
et 
$$
(1 \otimes c^*)(comp^*)(f)=\sum f_ic^*(\alpha^i(f)) + \sum e_lc^*(\beta^l(f)).
$$
En utilisant le second diagramme commutatif, on obtient les \'egalit\'es 
$$c^*(\beta^l(f_i)) =  \sum _{k} \beta^l(e_k)\beta^k(f_i)  +  \sum _{j}
\beta^l(f_j)\alpha^j(f_i)$$
qui prouvent la stabilit\'e par composition.
N'ayant pas de stabilit\'e par l'inversion $i$, consid\'erons l'id\'eal $
\mathcal{I}_p + i^*\mathcal{I}_p$. Par construction il est stable par
inversion et reste contenu dans l'id\'eal de l'identit\'e. Nous venons de prouver qu'une partie de cet id\'eal est stable par composition. Pour prouver que l'autre partie l'est aussi on introduit l'application $i\widehat{\times}i$ d\'efinie de $J^*(\Delta)\times_{\Delta}J^*(\Delta)$ dans $J^*(\Delta)\times_{\Delta}J^*(\Delta)$ par :
$$
\left((x,y,\ldots),(y,z,\ldots)\right) \to \left((z,y,\ldots),(y,x,\ldots)\right).  
$$
Cette fl\`eche induit un morphisme $(i\widehat{\times}i)^*$ v\'erifiant
$(i\widehat{\times}i)^*c^* =c^*i^*$. En l'appliquant \`a l'\'egalit\'e donnant la
stabilit\'e par composition de $\I_p$, on a :
$$
c^*(i^*(\beta^l(f_i))) =  \sum _{k}i^*(\beta^k(f_i))i^*(\beta^l(e_k))  +  \sum
_{j}i^*( \alpha^j(f_i)) i^*(\beta^l(f_j))
$$ 
ce qui prouve la stabilit\'e par composition de $ \mathcal{I}_p + i^*\mathcal{I}_p$.
%Le diagramme commutatif traduisant que l'application $i : \mathcal{O}_{J_k^*(\Delta)} \to \mathcal{O}_{J_k^*(\Delta)} $ est l'inverse pour la composition des jets donne $(diag^*)(1 \otimes i^*)(c^*)(\beta^l(f_i))=0$. Appliqu\'ee \`a l'\'egalit\'e obtenue pr\'ec\'edemment on obtient :
%$$
%0 =  \sum _{k}\beta^k(f_i) i^*(\beta^l(e_k)) +  \sum _{j} \alpha^j(f_i)i^*(\beta^l(f_j)).
%$$
Quitte \`a multiplier par les d\'enominateurs des $\beta^l(f_j)$, l'id\'eal
$\mathcal{I}_p+ i^*\mathcal{I}_p$ est inclus dans $\O_{J_p^*(\Delta)}$. Nous venons
de prouver qu'il d\'ecrit un sous-groupo\"ide de Lie de $J_p^*(\Delta)$,
la stabilit\'e par composition n'\'etant v\'erifi\'ee qu'en dehors du lieu des
z\'eros de ces d\'enominateurs. Le th\'eor\`eme de prolongement de
B. Malgrange \cite{malgrangegroupdegalois} assure que l'id\'eal r\'eduit qu'il  engendre diff\'erentiablement donne un $\mathcal{D}$-groupo\"ide de Lie.\\
V\'erifions que ce $\mathcal{D}$-groupo\"ide de Lie est admissible pour le
feuilletage. Les automorphismes locaux du feuilletage qui se factorisent en
l'identit\'e sur la transverse agissent comme l'identit\'e sur
$\mathcal{O}_{\F}$ donc laissent invariants tous ses id\'eaux
diff\'erentiels. Ils sont donc tous solutions du $\mathcal{D}$-groupo\"ide
de Lie que nous venons de construire.\\
V\'erifions enfin que ce \dgl\ est propre. On choisit des coordonn\'ees
transverses $t$ et des coordonn\'ees tangentes $z$.  Les \'equations
d'int\'egrales premi\`eres s'\'ecrivent localement $\frac{\partial H_i}{\partial
  z_j}=0$. Les \'equations suppl\'ementaires de l'id\'eal premier compatible
sont donc des \'equations que l'on peut supposer, quitte \`a effectuer les
substitutions n\'ecessaires, uniquement en les d\'eriv\'es des $H_i$ par rapport aux $t_k$. Cet id\'eal ne peut \^etre invariant par n'importe quelle transformation en $t$ : le $\mathcal{D}$-groupo\"ide de Lie construit est donc diff\'erent de celui de tous les automorphismes du feuilletage.
\end{preuve}

\begin{lemme}
Soit $\F$ un germe de feuilletage de codimension un dont le groupo\"ide de
Galois est propre. Il existe une int\'egrale premi\`ere de $\F$ dans une
extension fortement normale de $\M_{\Delta}$. 
\end{lemme}
\begin{preuve} En codimension un le rang transverse d'un \dgl\ admissible
  propre est fini.
La preuve de ce lemme se fait au cas par cas en discutant suivant le rang
transverse du groupo\"ide de Galois du feuilletage. On sait d'apr\`es le
th\'eoreme \ref{typedetransc} qu'il existe dans ces cas des int\'egrales
premi\`eres particuli\`eres. Ces int\'egrales premi\`eres vont nous donner
des id\'eaux diff\'erentiels premiers particuliers de $\O_\F$. Les corps des
fractions des quotients de $\O_\F$ par ces id\'eaux nous donneront dans chaque
cas une extension fortement normale contenant une int\'egrale premi\`ere.

\medskip
\textbf{Les feuilletages m\'eromorphiquement int\'egrables.}\rm \\
Lorsque le rang transverse est nul, le groupo\"ide de Galois est non transitif et
il existe une int\'egrale premi\`ere m\'eromorphe. L'extension est
$\M_{\Delta}$ et le morphisme est celui qui \`a $H$ associe une int\'egrale premi\`ere m\'eromorphe. 

\medskip
\indent \textbf{Les feuilletages Darboux-int\'egrables.}\rm \\
Lorsque le rang transverse est \'egal \`a un, rappelons comment on a
construit une int\'egrale premi\`ere de type Darboux. Pour $\Gamma \in Aut(\mathcal{F}_\omega)$, on note $f_\Gamma$ la fonction d\'efinie par $\Gamma^* \omega = f_\Gamma \omega$. On notera aussi $\omega=\sum w_idx_i$ dans des coordonn\'ees fix\'ees. Le groupo\"ide de Galois d'un feuilletage Darboux-int\'egrable est de la forme :
$$
\Gamma^*\omega \wedge \omega =0 \text{ et } m\circ\Gamma f_\Gamma^k=m
$$
pour une fonction m\'eromorphe $m$ et un entier $k$. Le syst\`eme d\'efini par les \'equations :
$$
(dH)^{\otimes k}-m(x)\omega^{\otimes k}=0
$$
fournit un id\'eal diff\'erentiel premier $\J$ de $\mathcal{O}_{\F}$. Le corps des
fractions $\mathcal{K}$ du quotient $\mathcal{O}_{\F}/\mathcal{J}$ est de degr\'e de transcendance un. Pour
prouver qu'il s'agit d'une extension fortement normale, prenons $\E$ une
extension diff\'erentielle de $\M_{\Delta}$ et $\sigma_1$ et $\sigma_2$ deux
plongements de $\mathcal{K}$ dans $\E$. Ces plongements sont compl\`etement
d\'etermin\'es par $H_1=\sigma_1 H$ et $H_2 = \sigma_2 H$, o\`u par abus de
notation $H$ d\'esigne aussi son image dans $\mathcal{K}$. Les \'el\'ements
$H_1$ et $H_2$ de $\E$ v\'erifient tous les deux l'\'equation engendrant
$\J$. Il existe donc une racine $k$-i\`eme de l'unit\'e $\theta$ et une
constante $b$ de $\E$ telles que $H_1=\theta H_2 +b$. Ayant obtenu une
expression rationnelle de $H_1$ en fonction de $H_2 $ \`a coefficients dans
les constantes de $\E$, on a $\sigma_1 \K \cdot \E^c=\sigma_2 \K \cdot \E^c$. Le corps $\mathcal{K}$ est donc une extension fortement normale de $\mathcal{M}_{\Delta}$ de groupe de Galois les transformations $x \to \theta x +b$ de la droite affine.
%est automorphe pour le groupo\"ide de Galois de $\mathcal{F}_\omega$. Les solutions de ce syst\`eme auront une monodromie de la forme $\varphi : H\to \theta H+b$ avec $\theta^k=1$. Le groupe de Drach du feuilletage est alors d\'efini par $\varphi'^k=1$. Pour obtenir la r\'esolvante de J. Drach, il suffit d\'ecrire le syst\`eme d'\'equations \`a r\'esoudre pour trouver $m$. \\
%Si il existe un entier $k$ et  une solution m\'eromorphe aux \'equations
%$$
%\frac{\partial m}{\partial x_j}w_i-\frac{\partial m}{\partial x_i}w_j=km\left(\frac{\partial w_i}{\partial x_j}-\frac{\partial w_j}{\partial x_i} \right) \text{ pour }1 \leq i,j \leq n
%$$   
%alors $m^{1/k}$ est un facteur int\'egrand $\omega$.

\medskip
\textbf{Les feuilletages Liouville-int\'egrables.}\rm \\
Lorsque le rang transverse de $\F_\omega$ est deux, il existe une forme ferm\'ee $\alpha$ v\'erifiant $d \omega = \omega \wedge \alpha $. 
Le groupo\"ide de Galois d'un tel feuilletage a pour \'equations :
$$
\Gamma^*\omega \wedge \omega =0 \text{ et } \Gamma^*\alpha + \frac{df_\Gamma}{f_\Gamma}=\alpha.
$$
On construit des int\'egrales premi\`eres particuli\`eres de ce feuilletage en r\'esolvant successivement 
$$
\frac{dF}{F}=\alpha \text{ puis } dH=F \omega.
$$
Dans des coordonn\'ees on \'ecrit $\omega= \sum w_i dx_i$ et $\alpha=\sum a_i
dx_i $. Le syst\`eme d'\'equations correspondant est :
$$
%\begin{array}{c}  
\frac{\frac{\partial H}{\partial x_i}}{w_i}=\frac{\frac{\partial H}{\partial x_j}}{w_j}\text{ ,  }\ \
\frac{\frac{\partial^2 H}{\partial x_i^2}}{\frac{\partial H}{\partial x_i}}=a_i + \frac{\frac{\partial w_i}{\partial x_i}}{w_i} \text{ pour }0\leq i,j \leq n.
%\end{array}
$$
Ce syst\`eme donne un id\'eal diff\'erentiel $\J$ de $\mathcal{O}_{\F}$. Le
corps des fractions du quotient $\mathcal{K}$ est de degr\'e de transcendance
deux. \'Etant donn\'ees deux solutions $H_1$ et $H_2$ des \'equations ci-dessus dans une extension de
$\mathcal{M}_{\Delta}$, on a $H_1= a H_2 +b$ avec $a$ et $b$ deux constantes de l'extension. Le corps $\mathcal{K}$ est donc une extension fortement normale de $\mathcal{M}_{\Delta}$ de groupe de Galois les transformations $x \to a x +b$ de la droite affine.
%Ce syst\`eme est bien automorphe car si $H$ est une solution, $\Gamma^*H$ v\'erifie $d(\Gamma^*H)=\Gamma^*F f_\Gamma \omega$ et $\frac{d(\Gamma^*F f_\Gamma)}{\Gamma^*F f_\Gamma}=\Gamma^*\alpha + \frac{df_\Gamma}{f_\Gamma}=\alpha$. De plus Si $H_1$ et $H_2$ sont deux solutions tels que $H_1=\Gamma^*H_2$, $\Gamma$ est solution du groupo\"ide de Galois du feuilletage. Les int\'egrales premi\`eres solutions de ce syst\`eme sont \`a monodromie affine. Le groupe de Drach est d\'efini par l'\'equation $\frac{\varphi''}{\varphi'}=0$. Pour trouver les r\'esolvantes de Drach, on commence par choisir une forme m\'eromorphe $\omega$ donnant le feuilletage v\'erifiant $w_1=1$, puis on \'ecrit $d\omega = \omega \wedge d\omega (\frac{\partial}{\partial x_1},\ )$. Ceci nous donne une forme $\alpha= \sum \frac{\partial w_i}{\partial x_1}dx_i$ v\'erifiant la premi\`ere \'equation de Godbillon-Vey. Il existe une suite de Godbillon-Vey de longueur deux si et seulement si il existe une fonction m\'eromorphe $R$ v\'erifiant $d(\alpha +R \omega )=0$. Lorsqu'on \'ecrit les \'equations correspondantes, on obtient la r\'esolvante :
%$$ 
%\frac{\partial^2 w_i}{\partial x_1 \partial x_j} - \frac{\partial^2 w_j}{\partial x_1 \partial x_i} + \frac{\partial(Rw_i)}{\partial x_j} - \frac{\partial(Rw_j)}{\partial x_i} \text{ pour }0\leq i,j \leq n.
%$$ 
%\\

\medskip
\textbf{Les feuilletages Riccati-int\'egrables.} \rm \\
Lorsque le rang transverse de $\F_\omega$ est trois, ce feuilletage admet une suite de Godbillon-Vey de longueur trois : $(\omega, \alpha, \beta)$. On construit des int\'egrales premi\`eres particuli\`eres en r\'esolvant la suite d'\'equations :
$$
\begin{array}{rl}
dG & =\frac{G^2}{2}\omega +G\alpha +\beta \\
\frac{dF}{F} & =G \omega +\alpha \\
dH & =F \omega.

\end{array} 
$$
Dans des coordonn\'ees on \'ecrit $\omega= \sum w_i dx_i$, $\alpha=\sum a_i
dx_i$, $\beta=\sum b_i dx_i$.
Le syst\`eme d'\'equations aux d\'eriv\'ees partielles correspondant est :

$$
\begin{array}{c}
\frac{\frac{\partial H}{\partial x_i}}{w_i}=\frac{\frac{\partial H}{\partial x_j}}{w_j} \\
\frac{1}{w_i}\left( \frac{\frac{\partial^2 H}{\partial x_i^2}}{\frac{\partial H}{\partial x_i}}-\frac{\frac{\partial w_i}{\partial x_i}}{w_i}-a_i \right)= \frac{1}{w_j}\left( \frac{\frac{\partial^2 H}{\partial x_j^2}}{\frac{\partial H}{\partial x_j}}-\frac{\frac{\partial w_j}{\partial x_j}}{w_j}-a_j \right)
\end{array} 
$$\label{sectiondiffeo}
\begin{multline}
\hspace*{1cm} \frac{\partial}{\partial x_i}\left(\frac{\frac{\partial^2 H}{\partial x_i^2}}{\frac{\partial H}{\partial x_i}}\right)-\frac{1}{2}\left(\frac{\frac{\partial^2 H}{\partial x_i^2}}{\frac{\partial H}{\partial x_i}}\right)^2= w_i\left(\left(\frac{\frac{\partial w_i}{\partial x_i}}{w_i}\right)-\frac{1}{2}\left(\frac{\frac{\partial w_i}{\partial x_i}}{w_i}\right)^2\right)  \nonumber \\ 
+\frac{\partial a_i}{\partial x_i}-\frac{1}{2}a_i^2 +b_i-a_i\frac{\frac{\partial w_i}{\partial x_i}}{w_i}. \hspace*{1cm} \nonumber
\end{multline}

%Le groupo\"ide de Galois d'un tel feuilletage a pour \'equations :
%$$
%\Gamma^*\omega \wedge \omega =0 \text{ et } f_\Gamma\Gamma^*\beta + dg_\Gamma = \beta + g_\gamma \alpha + \frac{1}{2}g^2_\Gamma \omega
%$$
Ces \'equations donnent un id\'eal diff\'erentiel $\J$ de
$\mathcal{O}_{\F}$. Le corps des fractions du quotient $\mathcal{K}$ est de
degr\'e de transcendance trois. \'Etant donn\'ees deux solutions $H_1$ et
$H_2$ des \'equations ci-dessus dans une extension de $\mathcal{M}_{\Delta}$ on a $H_1= \frac{a H_2 +b}{c H_2 +d}$ avec $a,b,c,d$ quatre constantes de l'extension. Le corps $\mathcal{K}$ est donc une extension fortement normale de $\mathcal{M}_{\Delta}$ de groupe de Galois les transformations homographiques de la droite projective.
\end{preuve}

%ùùùùùùùùùùùùùùùùùùùùùùùùùùùùùùùùùùùùùùùùùùùùùùùùùùùùùùùùùùùùùùùùùùùùùùùùùùùùùùùùùùùùùùùùùùùùùùùùùùùùùùùùùùùùùùùùùùùùùùùùùùùùùùùùùùùùùùùùùùùùùùùùùùùùùùùùùùùùùùùùùùùùùùùùùùùùùùùùùùùùùùùùùùùùùùùùùùùùùùùùùùùùùùùùùùùùùùùùùùùùùùùùùùùùùùùùùùùùùùùùùùùùùùùùùùùùùùùùùùùùùùùùùùùùùùùùùùùùùùùùùùùùùùùùùùùùùùùùùùùùùùùùùùùùùùùùùùùùùùùùùùùùùùùùùùùùùùùùùùùùùùùùùùùùùùùùùùùùùùùùùùùùùùùùùùùùùùùùùùùùùùùùùùùùùùùùùùùùù
%ù          LA BIBLIO
%ù
%ùùùùùùùùùùùùùùùùùùùùùùùùùùùùùùùùùùùùùùùùùùùùùùùùùùùùùùùùùùùùùùùùùùùùùùùùùùùùùùùùùùùùùùùùùùùùùùùùùùùùùùùùùùùùùùùùùùùùùùùùùùùùùùùùùùùùùùùùùùùùùùùùùùùùùùùùùùùùùùùùùùùùùùùùùùùùùùùùùùùùùùùùùùùùùùùùùùùùùùùùùùùùùùùùùùùùùùùùùùùùùùùùùùùùùùùùùùùùùùùùùùùùùùùùùùùùùùùùùùùùùùùùùùùùùùùùùùùùùùùùùùùùùùùùùùùùùùùùùù

\end{document}